\documentclass[12pt]{article}

\usepackage[T1]{fontenc}
\usepackage{amsmath}
\usepackage{amssymb}
\usepackage{latexsym}
\usepackage{textcomp}
\usepackage{diagrams}
\diagramstyle{h=25pt,w=28pt,tight,scriptlabels,midshaft,PostScript=dvips,nohug,shortfall=4pt,heads=littlevee}
\def\rhcvee{\mkern-10mu\hbox{\footnotesize \raisebox{0.5pt}{$\succ$}}}%%
\def\lhcvee{\hbox{\footnotesize \raisebox{0.5pt}{$\prec$}}\mkern-10mu}%%
\def\dhcvee{\vboxtoz{\vss\hbox{\footnotesize $\curlyvee$}\kern0pt}}%%
\def\uhcvee{\vboxtoz{\hbox{\footnotesize $\curlywedge$}\vss}}%%
\newarrowhead{littlevee}\rhcvee\lhcvee\dhcvee\uhcvee

\newarrow{Lig}{=}{=}{=}{=}{=}

\usepackage[ps,dvips,all]{xy}
\xyoption{ps}
\SelectTips{cm}{}
\usepackage[noBBpl]{mathpazo}
\renewcommand{\texttt}[1]{{\fontfamily{pcr}\fontseries{m}\fontshape{n}%
\selectfont #1}}

\providecommand{\overskrift}[1]{\par\noindent\relax{\LARGE #1}\par\bigskip}

\newcommand{\hovedfont}{\normalfont\bfseries}
\usepackage{ntheorem}
\theoremheaderfont{\hovedfont}

\theoremseparator{.}
    \theorembodyfont{\normalfont\itshape}
    \theoremstyle{plain}
\newtheorem{satz}{Theorem}

    \theoremstyle{change}
\newtheorem{lemma}{Lemma}[section]
\newtheorem{prop}[lemma]{Proposition}
\newtheorem{cor}[lemma]{Corollary}
    \theorembodyfont{\normalfont}

\newtheorem{BM}[lemma]{Remark}

\newtheorem{taller}[lemma]{$\!\!$}
\newenvironment{blanko}[1]%
{\begin{taller}{\hovedfont #1}\normalfont}%
{\end{taller}}
{%
\begin{list}{\em Definition. }%
{\setlength{\labelsep}{0mm}\setlength{\leftmargin}{0mm}%
\setlength{\labelwidth}{0mm}\setlength{\listparindent}{\parindent}%
\setlength{\parsep}{\parskip}\setlength{\partopsep}{0mm}}%
\item%
}%
{%
\end{list}%
}
\newenvironment{dem}%
{%
\begin{list}{\em Proof. }%
{\setlength{\labelsep}{0mm}\setlength{\leftmargin}{0mm}%
\setlength{\labelwidth}{0mm}\setlength{\listparindent}{\parindent}%
\setlength{\parsep}{\parskip}\setlength{\partopsep}{0mm}}%
\item%
}%
{%
\qed\end{list}%
}
\newenvironment{dem*}[1]%
{%
\begin{list}{\em #1 }%
{\setlength{\labelsep}{0mm}\setlength{\leftmargin}{0mm}%
\setlength{\labelwidth}{0mm}\setlength{\listparindent}{\parindent}%
\setlength{\parsep}{\parskip}\setlength{\partopsep}{0mm}}%
\item%
}%
{%
\qed\end{list}%
}
{%
\begin{list}{\em Proof. }%
{\setlength{\labelsep}{0mm}\setlength{\leftmargin}{0mm}%
\setlength{\labelwidth}{0mm}\setlength{\listparindent}{\parindent}%
\setlength{\parsep}{\parskip}\setlength{\partopsep}{0mm}}%
\item%
}%
{%
\qed\end{list}%
}
\newenvironment{bevis*}[1]%
{%
\begin{list}{\em #1 }%
{\setlength{\labelsep}{0mm}\setlength{\leftmargin}{0mm}%
\setlength{\labelwidth}{0mm}\setlength{\listparindent}{\parindent}%
\setlength{\parsep}{\parskip}\setlength{\partopsep}{0mm}}%
\item%
}%
{%
\qed\end{list}%
}

\newenvironment{blanko*}[1]%
{%
\begin{list}{\bf {#1} }%
{\setlength{\labelsep}{0mm}\setlength{\leftmargin}{0mm}%
\setlength{\labelwidth}{0mm}\setlength{\listparindent}{\parindent}%
\setlength{\parsep}{\parskip}\setlength{\partopsep}{0mm}}%
\item%
}%
{%
\end{list}%
}

\newcounter{dummycounter}

\newenvironment{punkt-i}%
{%
	\begin{list}%
	{(\roman{dummycounter})}%
	{\usecounter{dummycounter}%
	\setlength{\itemsep}{0em}\setlength{\parsep}{0em}\setlength{\topsep}{0em}%
	\setlength{\itemindent}{0em}\setlength{\labelwidth}{1.8em}%
	\setlength{\labelsep}{0.6em}\setlength{\leftmargin}{2.4em}}%
}%
{\end{list}}

\usepackage{mathrsfs}
\newcommand{\CC}{\mathscr{C}}
\newcommand{\DD}{\mathscr{D}}
\newcommand{\EE}{\mathscr{E}}
\newcommand{\GG}{\mathscr{G}}
\newcommand{\UU}{\mathscr{U}}

\providecommand{\lastUpdate}[1]{#1}

\newcommand{\Nat}{\operatorname{Nat}}
\newcommand{\id}{\operatorname{id}}
\newcommand{\Id}{\operatorname{Id}}

\newcommand{\comm}{\copyright}

\newcommand{\tensor}{\otimes}
\newcommand{\isopil}{\stackrel{\raisebox{0.1ex}[0ex][0ex]{\(\sim\)}}%
			{\raisebox{-0.15ex}[0.28ex]{\(\rightarrow\)}}}
\newcommand{\Hom}{\operatorname{Hom}}
\providecommand{\qed}{\hspace*{\fill}\nolinebreak[1]\hspace*{\fill}$\Box$}
\newcommand{\df}{\: {\raisebox{0.255ex}{\normalfont\scriptsize :\!\!}}=}

\newcommand{\cel}[1]{\ensuremath{\mathsf{#1}}}
\newcommand{\Aa}{\text{\rm \textsf{\AA}}}

\newcommand{\inv}{^{\scriptscriptstyle -1}}

\newcommand{\loft}{^{\text{\tiny \rm{left}}}}
\newcommand{\rught}{^{\text{\tiny \rm{right}}}}

\newcommand{\CCa}{\CC^{\mathbf{2}}}

\newcommand{\impl}{$\Rightarrow$}

\def\vspec#1{\special{ps:#1}}%  passes #1 verbatim to the output
\def\rotstart#1{\vspec{gsave currentpoint currentpoint translate
	#1 neg exch neg exch translate}}% #1 can be any origin-fixing transformation
\def\rotfinish{\vspec{currentpoint grestore moveto}}% gets back in synch

\def\psr#1{\rotstart{36.87 rotate}\hbox to0pt {\vsize 0pt \hss\(\scriptstyle #1\)\hss}\rotfinish}
\def\sdf#1{\rotstart{-36.87 rotate}\hbox to0pt {\vsize 0pt \hss\(\scriptstyle #1\)\hss}\rotfinish}

\def\scaleFactor{10} % MUST BE AN INTEGER

\newcommand{\dropglob}[1]{%
\setlength{\unitlength}{0.003\DiagramCellWidth}
\multiply \unitlength by \scaleFactor
\begin{picture}(0,0)(0,0)
\qbezier(-28,-4)(0,-18)(28,-4)
\put(0,-14){\makebox(0,0)[t]{$\scriptstyle {#1}$}}
\put(28.6,-3.7){\vector(2,1){0}}
\end{picture}
}

\newcommand{\topglob}[1]{%
\setlength{\unitlength}{0.003\DiagramCellWidth}
\multiply \unitlength by \scaleFactor
\begin{picture}(0,0)(0,0)
\qbezier(-28,11)(0,25)(28,11)
\put(0,21){\makebox(0,0)[b]{$\scriptstyle {#1}$}}
\put(28.6,10.7){\vector(2,-1){0}}
\end{picture}
}

\newcommand{\leftglob}[1]{%
\setlength{\unitlength}{0.003\DiagramCellHeight}
\multiply \unitlength by \scaleFactor
\begin{picture}(0,0)(0,-4)
\qbezier(-9,-25)(-23,0)(-9,25)
\put(-20,0){\makebox(0,0)[r]{$\scriptstyle {#1}$}}
\put(-8.7,-25.6){\vector(1,-2){0}}
\end{picture}
}
\newcommand{\rightglob}[1]{%
\setlength{\unitlength}{0.003\DiagramCellHeight}
\multiply \unitlength by \scaleFactor
\begin{picture}(0,0)(0,-4)
\qbezier(9,-25)(23,0)(9,25)
\put(20,0){\makebox(0,0)[l]{$\scriptstyle {#1}$}}
\put(8.7,-25.6){\vector(-1,-2){0}}
\end{picture}
}

\newcommand{\thindropglob}[1]{%
\setlength{\unitlength}{0.003\DiagramCellWidth}
\multiply \unitlength by \scaleFactor
\begin{picture}(0,0)(0,0)
\qbezier(-26,-1)(0,-9)(26,-1)
\put(0,-7){\makebox(0,0)[t]{$\scriptstyle {#1}$}}
\put(26.6,-0.7){\vector(3,1){0}}
\end{picture}
}

\newcommand{\thintopglob}[1]{%
\setlength{\unitlength}{0.003\DiagramCellWidth}
\multiply \unitlength by \scaleFactor
\begin{picture}(0,0)(0,0)
\qbezier(-26,6)(0,14)(26,6)
\put(0,12){\makebox(0,0)[b]{$\scriptstyle {#1}$}}
\put(26.6,5.7){\vector(3,-1){0}}
\end{picture}
}

\newcommand{\thinleftglob}[1]{%
\setlength{\unitlength}{0.003\DiagramCellHeight}
\multiply \unitlength by \scaleFactor
\begin{picture}(0,0)(0,-4)
\qbezier(-3,-25)(-16,0)(-3,25)
\put(-13,0){\makebox(0,0)[r]{$\scriptstyle {#1}$}}
\put(-2.5,-25.6){\vector(1,-2){0}}
\end{picture}
}

\newcommand{\thinrightglob}[1]{%
\setlength{\unitlength}{0.003\DiagramCellHeight}
\multiply \unitlength by \scaleFactor
\begin{picture}(0,0)(0,-4)
\qbezier(3,-25)(16,0)(3,25)
\put(13,0){\makebox(0,0)[l]{$\scriptstyle {#1}$}}
\put(2.7,-25.6){\vector(-1,-2){0}}
\end{picture}
}

\newcommand{\diagonaldropglob}[1]{%
\setlength{\unitlength}{0.003\DiagramCellWidth}
\multiply \unitlength by \scaleFactor
\psr{\begin{picture}(0,0)(0,0)
\qbezier(-34,1)(0,-8)(30,1)
\put(0,-6){\makebox(0,0)[t]{$\scriptstyle {#1}$}}
\put(30.6,1.3){\vector(3,1){0}}
\end{picture}}
}

\newcommand{\diagonaltopglob}[1]{%
\setlength{\unitlength}{0.003\DiagramCellWidth}
\multiply \unitlength by \scaleFactor
\psr{\begin{picture}(0,0)(0,0)
\qbezier(-34,6)(0,16)(30,6)
\put(0,13){\makebox(0,0)[b]{$\scriptstyle {#1}$}}
\put(30.6,5.7){\vector(3,-1){0}}
\end{picture}}
}

\newcommand{\bigdropglob}[1]{%
\setlength{\unitlength}{0.003\DiagramCellWidth}
\multiply \unitlength by \scaleFactor
\begin{picture}(0,0)(0,0)
\qbezier(-28,-4)(0,-26)(28,-4)
\put(0,-17){\makebox(0,0)[t]{$\scriptstyle {#1}$}}
\put(28.6,-3.7){\vector(1,1){0}}
\end{picture}
}

\newcommand{\bigtopglob}[1]{%
\setlength{\unitlength}{0.003\DiagramCellWidth}
\multiply \unitlength by \scaleFactor
\begin{picture}(0,0)(0,0)
\qbezier(-28,11)(0,33)(28,11)
\put(0,24){\makebox(0,0)[b]{$\scriptstyle {#1}$}}
\put(28.6,10.7){\vector(1,-1){0}}
\end{picture}
}

\newcommand{\straighthor}[1]{%
\setlength{\unitlength}{0.003\DiagramCellWidth}
\multiply \unitlength by \scaleFactor
\begin{picture}(0,0)(0,0)
\put(15,4){\makebox(0,0)[b]{$\scriptstyle {#1}$}}
\put(-21,3){\vector(1,0){44}}
\end{picture}
}

\newcommand{\lift}[2]{%
\setlength{\unitlength}{1pt}
\begin{picture}(0,0)(0,0)
\put(0,{#1}){\makebox(0,0)[b]{${#2}$}}
\end{picture}
}

\newcommand{\shuft}[2]{%
\setlength{\unitlength}{1pt}
\begin{picture}(0,0)(0,0)
\put({#1},0){\makebox(0,0)[b]{${#2}$}}
\end{picture}
}

\newcommand{\NEdiagonaltopglob}[1]{%
\setlength{\unitlength}{0.003\DiagramCellWidth}
\multiply \unitlength by \scaleFactor
\sdf{\begin{picture}(0,0)(0,0)
\qbezier(-32,7)(2,19)(36,7)
\put(2,16){\makebox(0,0)[b]{$\scriptstyle {#1}$}}
\put(36.6,7){\vector(3,-1){0}}
\end{picture}}
}

\newcommand{\NEdiagonaldropglob}[1]{%
\setlength{\unitlength}{0.003\DiagramCellWidth}
\multiply \unitlength by \scaleFactor
\sdf{\begin{picture}(0,0)(0,0)
\qbezier(36,0)(2,-12)(-32,0)
\put(2,-9){\makebox(0,0)[t]{$\scriptstyle {#1}$}}
\put(36.6,0.6){\vector(3,1){0}}
\end{picture}}
}

\newcommand{\NESWlabel}[1]{%
\setlength{\unitlength}{0.003\DiagramCellWidth}
\multiply \unitlength by \scaleFactor
\sdf{\begin{picture}(0,0)(0,0)
\put(2,8){\makebox(0,0)[t]{$#1$}}
\end{picture}}
}
\newcommand{\NWSElabel}[1]{%
\setlength{\unitlength}{0.003\DiagramCellWidth}
\multiply \unitlength by \scaleFactor
\psr{\begin{picture}(0,0)(0,0)
\put(2,8){\makebox(0,0)[t]{$#1$}}
\end{picture}}
}
\makeatletter

\renewcommand{\tableofcontents}{%
   \begin{center}
\begin{minipage}{110mm}
   \begin{center}
		\bf{\contentsname}
	\end{center}
   \footnotesize
   \begin{center}
		\@starttoc{toc}
	\end{center}	
\end{minipage}
	\end{center}
	\addvspace{3em \@plus\p@}
}

\renewcommand{\section}{\@startsection {section}{1}{\z@}%
{-3.5ex \@plus -1ex \@minus -.2ex}%
{2.3ex \@plus.2ex}%
{\normalfont\large\bfseries}}

\renewcommand*{\l@section}[2]{%
  \ifnum \c@tocdepth >\z@
    \addpenalty\@secpenalty
    \setlength\@tempdima{1.5em}%
    \begingroup
      \parindent \z@ \rightskip \@pnumwidth
      \parfillskip -\@pnumwidth
      \leavevmode %\bfseries
      \advance\leftskip\@tempdima
      \hskip -\leftskip
      #1\nobreak\hfil \nobreak\hb@xt@\@pnumwidth{\hss #2}\par
    \endgroup
  \fi}

\makeatother

\makeatletter
\renewcommand{\ps@headings}
	{\setlength{\headheight}{13pt}%
	 \setlength{\headsep}{12pt}%
	 \renewcommand{\@oddhead}{\parbox{\textwidth}{%
			\small
			\texttt{\jobname.tex \ \ \ \lastUpdate {2009-07-18 09:24}
			\hfill [\thepage/\pageref{lastpage}]}
			\\ \rule[8pt]{\textwidth}{0.3pt}}%
	 }
	\renewcommand{\@oddfoot}{}
	\renewcommand{\@evenfoot}{}%
}
\makeatother

\begin{document}

\pagestyle{headings}

\vspace*{36pt}

\overskrift{Coherence for weak units}

\vspace*{8pt}

\noindent
{\large \textsc{Andr\'e Joyal} and \textsc{Joachim Kock}}

\vspace{4pt}

\noindent
{\small Universit\'e du Qu\'ebec \`a Montr\'eal}

% Keywords: monoidal $2$-categories, units, coherence.

% MSC: 18D05 Double categories, $2$-categories, bicategories and generalizations

% Comments: LaTeX 37 pages; does not compile with pdflatex due to some ps 
% rotations.  Minor typographical shortcomings.

\vspace{16pt}

\hrule

\vspace{8pt}

\footnotesize

\noindent \textbf{Abstract.} We define weak units in a semi-monoidal
$2$-category $\CC$ as cancellable pseudo-idempotents: they are pairs
$(I,\alpha)$ where $I$ is an object such that tensoring with $I$ from either
side constitutes a biequivalence of $\CC$, and $\alpha: I \tensor I \to I$ is an
equivalence in $\CC$.  We show that this notion of weak unit has coherence built
in: Theorem~\ref{thmA}: $\alpha$ has a canonical associator $2$-cell, which
automatically satisfies the pentagon equation.  Theorem~\ref{thmB}: every
morphism of weak units is automatically compatible with those associators.
Theorem~\ref{thmC}: the $2$-category of weak units is contractible if non-empty.
Finally we show (Theorem~\ref{thmE}) that the notion of weak unit is equivalent
to the notion obtained from the definition of tricategory: $\alpha$ alone
induces the whole family of left and right maps (indexed by the objects),
as well as the whole family of Kelly $2$-cells (one for
each pair of objects), satisfying the relevant coherence axioms.

\normalsize

\vspace{8pt}

\hrule

\bigskip
  
% \tableofcontents

%%%%%%%%%%%%%%%%%%%%%%%%%%%%%%%%%%%%%%%%%%%%%%%%%%
\section*{Introduction}
%%%%%%%%%%%%%%%%%%%%%%%%%%%%%%%%%%%%%%%%%%%%%%%%%%

The notion of tricategory, introduced by Gordon, Power, and
Street~\cite{Gordon-Power-Street} in 1995, seems still to represent the
highest-dimensional explicit weak categorical structure that can be manipulated
by hand (i.e.~without methods of homotopy theory), and is therefore an
important test bed for higher-categorical ideas.  In this work we investigate
the nature of weak units at this level.  While coherence for weak associativity
is rather well understood, thanks to the geometrical insight provided by the
Stasheff associahedra~\cite{Stasheff:1963}, coherence for unit structures is
more mysterious, and so far there seems to be no clear geometric pattern for the
coherence laws for units in higher dimensions.  Specific interest in weak units
stems from Simpson's conjecture~\cite{Simpson:9810}, according to which strict
$n$-groupoids with weak units should model all homotopy $n$-types.

In the present paper, working in the setting of a strict $2$-category $\CC$ with a
strict tensor product, we define a notion of weak unit by simple axioms that
involve only the notion of equivalence, and hence in principle make sense in all
dimensions.  Briefly, a weak unit is a cancellable pseudo-idempotent.  We
work out the basic theory of such units, and compare with the notion extracted
from the definition of tricategory.  In the companion paper {\em Weak units and
homotopy $3$-types}~\cite{Joyal-Kock:traintracks} we employ this notion of unit
to prove a version of Simpson's conjecture for $1$-connected homotopy $3$-types,
which is the first nontrivial case.  The strictness assumptions of the present
paper should be justified by that result.

\bigskip

By cancellable pseudo-idempotent we mean a pair $(I,\alpha)$ where $I$ is an
object in $\CC$ such that tensoring with $I$ from either side is an equivalence of
$2$-categories, and $\alpha: I \tensor I \isopil I$ is an equi-arrow (i.e.~an
arrow admitting a pseudo-inverse).  The remarkable fact about this definition is
that $\alpha$, viewed as a multiplication map, comes with canonical higher order
data built in: it possesses a canonical associator $\cel A$ 
which automatically satisfies the pentagon equation.
This is our Theorem~\ref{thmA}.  The point is that the arrow $\alpha$ alone,
thanks to the cancellability of $I$, induces all the usual structure of left and right
constraints with all the $2$-cell data that goes into them and the axioms they
must satisfy.  
% While all this left and right data is non-canonical,
% the associator $\cel A$ is canonical.

\bigskip

As a warm-up to the various constructions and ideas, we start out in
Section~\ref{sec:dim1} by briefly running through the corresponding theory for
cancellable-idempotent units in monoidal $1$-categories.  This theory has been
treated in detail in \cite{Kock:0507349}.
  
The rest of the paper is dedicated to the case of monoidal 
$2$-categories.  In Section~\ref{sec:dim2-main} we give the definitions
and state the main results:
Theorem~\ref{thmA} says that there is a canonical associator $2$-cell for
$\alpha$, and that this $2$-cell automatically satisfies the
pentagon equation.
Theorem~\ref{thmB} states  that unit morphisms automatically
are compatible with the associators of Theorem~\ref{thmA}.
Theorem~\ref{thmC} states that the $2$-category of units is
contractible if non-empty.  Hence, `being unital' is, up to homotopy, a property
rather than a structure. 

Next follow three sections dedicated to proofs of each of these three theorems.
In Section~\ref{sec:left-right} we show how the map $\alpha:II\isopil I$ alone
induces left and right constraints, which in turn are used to construct the
associator and establish the pentagon equation.  The left and right
constraints are not canonical, but surprisingly the associator does not depend
on the choice of them.  In Section~\ref{sec:catArr} we prove Theorem~\ref{thmB} by
interpreting it as a statement about units in the $2$-category of arrows, where
it is possible to derive it from Theorem~\ref{thmA}. In Section~\ref{sec:contractibility}
we prove Theorem~\ref{thmC}. The key ingredient is to use the left and right constraints
to link up all the units, and to show that the unit morphisms are precisely
those compatible with the left and right constraints; this makes them
`essentially unique' in the required sense.

In Section~\ref{sec:classical} we go through the basic theory of classical units
(i.e.~as extracted from the definition of
tricategory~\cite{Gordon-Power-Street}).  Finally, in
Section~\ref{sec:comparison} we show that the two notions of unit are
equivalent.  This is our Theorem~\ref{thmE}. A curiosity implied by the arguments in this
section is that the left and right axioms for the $2$-cell data in the 
Gordon-Power-Street
definition (denoted TA2 and TA3 in \cite{Gordon-Power-Street}) imply each other.
% an observation that may not have been made before.

(We have no Theorem D.)

\bigskip

This notion of weak units as cancellable idempotents is precisely what can be
extracted from the more abstract, Tamsamani-style, theory of fair
$n$-categories~\cite{Kock:0507116} by making an arbitrary choice of a fixed weak
unit.  In the theory of fair categories, the key object is a contractible space
of all weak units, rather than any particular point in that space, and handling
this space as a whole bypasses coherence issues.  However, for the sake of
understanding what the theory entails, and for the sake of concrete
computations, it is interesting to make a choice and study the ensuing coherence
issues, as we do in this paper.  The resulting approach is very much in the
spirit of the classical theory of monoidal categories, bicategories, and
tricategories, and provides some new insight to these theories.  To stress this
fact we have chosen to formulate everything from scratch in such classical
terms, without reference to the theory of fair categories.

In the case of monoidal $1$-categories, the cancellable-idempotent viewpoint on
units goes back to Saavedra~\cite{Saavedra}.  The importance of this viewpoint
in higher categories was first suggested by Simpson~\cite{Simpson:9810}, in
connection with his weak-unit conjecture.  He gave an ad hoc definition in this
style, as a mere indication of what needed to be done, and raised the question
of whether higher homotopical data would have to be specified.  The surprising
answer is, at least here in dimension $3$, that specifying $\alpha$ is enough,
then the higher homotopical data is automatically built in.

\bigskip

This paper was essentially written in 2004, in parallel with
\cite{Joyal-Kock:traintracks}.  We are ourselves to blame for the delay of
getting it out of the door.  The present form of the paper represents only half
of what was originally planned to go into the paper.  The second half should
contain an analysis of strong monoidal functors (along the lines of what was
meanwhile treated just in the $1$-dimensional case \cite{Kock:0507349}), and
also a construction of the `universal unit', hinted at in \cite{Kock:0507116}.
We regret that these ambitions should hold back the present material for so
long, and have finally decided to make this first part available {\em as is}, in
the belief that it is already of some interest and can well stand alone.

\begin{blanko*}{Acknowledgements.}
  We thank Georges Maltsiniotis for pointing out to us that the
  cancellable-idempotent notion of unit in dimension $1$ goes back to
  Saavedra~\cite{Saavedra}, and we thank Josep Elgueta for catching an error in
  an earlier version of our comparison with tricategories.  The first-named
  author was supported by the NSERC. The second-named author was very happy to
  be a CIRGET postdoc at the UQAM in 2004, and currently holds support from
  grants MTM2006-11391 and MTM2007-63277 of Spain.
\end{blanko*}

%%%%%%%%%%%%%%%%%%%%%%%%%%%%%%%%%%%%%%%%%%%%%%%%%%
\section{Units in monoidal categories}
%%%%%%%%%%%%%%%%%%%%%%%%%%%%%%%%%%%%%%%%%%%%%%%%%%
\label{sec:dim1}

It is helpful first briefly to recall the relevant results for
monoidal categories, referring the reader to \cite{Kock:0507349} for
further details of this case.

\begin{blanko}{Semi-monoidal categories.} 
%   (B\'enabou~\cite{Benabou:CR1963}, Mac
%   Lane~\cite{MacLane:naturalAssociativity}.)  
  A {\em semi-monoidal category} is a category $\CC$ equipped with a tensor
  product (which we denote by plain juxtaposition), i.e.~an associative functor
  \begin{eqnarray*}
    \CC\times\CC & \longrightarrow & \CC  \\
    (X,Y) & \longmapsto & XY .
  \end{eqnarray*}
  For simplicity we assume strict associativity,
  $X(YZ)=(XY)Z$.
%   This is really no loss of generality: all the arguments carry
%   over to the case of non-strict associativity --- just insert associators where
%   needed.
\end{blanko}

\begin{blanko}{Monoidal categories.}
  (Mac Lane~\cite{MacLane:naturalAssociativity}.)  A semi-monoidal category $\CC$ is a
  {\em monoidal category} when it is furthermore equipped with a distinguished
  object $I$ and natural isomorphisms
\begin{diagram}[w=4ex,tight]
I X & \rTo^{\lambda_X}  & X 
  &  \lTo^{\rho_X}  &  X I 
\end{diagram}
obeying the following rules (cf.~\cite{MacLane:naturalAssociativity}):
\begin{align} 
  % Saavedra (1.3.1.1)
  % Mac Lane (5.6)
  % Kelly (4)
  \lambda_I &= \rho_I \\
% \end{equation}
% \begin{equation} 
  % Saavedra (2.2.1.1 (b))
  % Mac Lane (5.3 (i))
  % Kelly (5)
  \lambda_{XY} &= \lambda_X Y \\
% \end{equation}
% \begin{equation} 
  % Saavedra (2.2.1.1 (a))
  % Mac Lane (5.7 (ii))
  % Kelly (7)
 \rho_{XY} &= X\rho_Y \\
% \end{equation}
% \begin{equation}
  % Saavedra (2.2.1.1 (c))
  % Mac Lane (5.7 (i))
  % Kelly (6)
  X\lambda_Y &= \rho_X Y \label{Kelly} 
\end{align}
% It was observed by Kelly~\cite{Kelly:MacLanesCoherence} that
% Axiom~(\ref{Kelly}) implies the other three axioms.
\end{blanko}
Naturality of $\lambda$ and $\rho$ implies
\begin{equation}\label{naturality1}
  \lambda_{IX} = I\lambda_X ,  \qquad \rho_{XI} = \rho_X I   ,
\end{equation}
independently of Axioms (1)--(4).

\begin{BM}
  Tensoring with $I$ from either side is an equivalence of categories.
\end{BM}

\begin{lemma}
  {\rm (Kelly~\cite{Kelly:MacLanesCoherence}.)}  Axiom~(4) implies axioms
  (1), (2), and (3).
\end{lemma}
\begin{dem}
    (4) implies (2):
    Since tensoring with $I$ on the left is an equivalence, it is 
    enough to prove $I \lambda_{XY} = I\lambda_X Y$.  But this
    follows from Axiom~(4) applied twice (swap $\lambda$ out 
    for a $\rho$ and swap back again only on the nearest factor):
    $$
    I\lambda_{XY} = \rho_I XY = I \lambda_X Y .
    $$
    Similarly for $\rho$, establishing (3).
    
    (4) and (2) implies (1):
  Since tensoring with $I$ on the right is an equivalence, it is 
  enough to prove $\lambda_I I = \rho_I I$.  But this follows from
  (2), (5), and (4):
  $$
  \lambda_I I = \lambda_{II} = I \lambda_I = \rho_I I .
  $$
  
  \vspace{-2em}
\end{dem}

% It was observed by Saavedra~\cite{Saavedra} that conversely (1), (2),
% and (3) imply (4).

The following alternative notion of unit object goes back to
Saavedra~\cite{Saavedra}.  A thorough treatment of the notion was given in
\cite{Kock:0507349}.

\begin{blanko}{Units as cancellable pseudo-idempotents.}\label{unit1}
  An object $I$ in a semi-monoidal category $\CC$ is called {\em cancellable} if the two
  functors $\CC\to\CC$
  \begin{eqnarray*}
    X & \longmapsto & IX\\ 
    X & \longmapsto & XI 
  \end{eqnarray*}
  are fully faithful.  
% REMARK: for the sake of defining units, in fact it is enough to
% be fully faithful with respect to isomorphisms!
% 
% (Saavedra calls such objects $1$-regular.)
  By definition, a {\em pseudo-idempotent} is
  an object $I$ equipped with an isomorphism $\alpha : II \isopil I$.
  Finally we define a {\em unit object} in $\CC$ to be a cancellable pseudo-idempotent.
\end{blanko}
% (This is called `reduced unit' by Saavedra~\cite{Saavedra}, Ch.~I, 1.3.2.)

\begin{lemma}\label{Kelly-lemma}
  {\rm \cite{Kock:0507349}}
  Given a unit object $(I,\alpha)$ in a semi-monoidal category $\CC$, for each
  object $X$ there are unique arrows
  \begin{diagram}[w=4ex,tight]
    I X & \rTo^{\lambda_X}  & X 
    &  \lTo^{\rho_X}  &  X I 
  \end{diagram}
  such that 
  \begin{eqnarray*}
    % Saavedra, (2.2.3.3)
    (\cel{L}) &  & I\lambda_X = \alpha X  \\
    % Saavedra, (2.2.3.2)
    (\cel{R}) &  & \rho_X I = X \alpha .
  \end{eqnarray*}
  The $\lambda_X$ and $\rho_X$ are isomorphisms and natural in $X$.
\end{lemma}

\begin{dem}
  Let $\mathbb{L}: \CC\to \CC$ denote the functor defined by tensoring with $I$
  on the left.  Since $\mathbb{L}$ is fully faithful, we have a bijection
  $$
  \Hom(IX,X) \to \Hom(IIX,IX).
  $$  
  Now take $\lambda_X$ to be the inverse image of $\alpha X$; it is an
  isomorphism since $\alpha X$ is.  Naturality follows by considering more
  generally the bijection
  $$
  \Nat(\mathbb{L},\id_\CC) \to \Nat(\mathbb{L}\circ \mathbb{L}, \mathbb{L}) 
  ;
  $$
  let $\lambda$ be the inverse image of the natural
  transformation whose components are $\alpha X$.
  Similarly on the right.
\end{dem}
% % Observation:  naturality implies that tensoring with $I$ on the left
% % (or on the right) is an equifunctor $\CC\to\CC$.  Indeed, the natural 
% % transformation $\lambda$ witnesses that tensoring with $I$ is
% % isomorphic to the identity functor.  Conversely, it is shown that
% % if left-tensoring with $I$ is an equifunctor, then every family of
% % $\lambda_X$ is natural:
% % Let $\lambda_X$
% % and $\lambda_Y$ be any two isomorphisms, and let $f:X\to Y$.
% % Then certainly we have a commutative square
% % \begin{diagram}[w=6ex,h=4.5ex,tight]
% % IIX & \rTo^{\alpha X}  & IX  \\
% % \dTo<{IIf}  &    & \dTo>{If}  \\
% % IIY  & \rTo_{\alpha Y}  & IY
% % \end{diagram}
% % Using the previous observation to substitute $\alpha X = I \lambda_X$
% % and $\alpha Y = I \lambda_Y$.  Hence
% % \begin{diagram}[w=6ex,h=4.5ex,tight]
% % IIX & \rTo^{I \lambda_X}  & IX  \\
% % \dTo<{IIf}  &    & \dTo>{If}  \\
% % IIY  & \rTo_{I \lambda_Y}  & IY
% % \end{diagram}
% % Now cancel the left-hand $I$ factor --- we can do this since the
% % functor `tensor with $I$' is an isomorphism of categories. 

\begin{lemma}\label{Kelly-dim1}
  {\rm \cite{Kock:0507349}}
  For $\lambda$ and $\rho$ as above, the Kelly axiom~\eqref{Kelly} holds:
  $$
  X \lambda_Y = \rho_X Y .
  $$
  Therefore, by Lemma~\ref{Kelly-lemma} a semi-monoidal category with
  a unit object is a monoidal category in the classical sense.
\end{lemma}

\begin{dem}
  In the commutative square
  \begin{diagram}[w=6ex,h=4.5ex,tight]
    XIIY & \rTo^{XI\lambda_Y}  & XIY  \\
    \dTo<{\rho_X I Y}  &    & \dTo>{\rho_X Y}  \\
    XIY  & \rTo_{X\lambda_Y}  & XY
  \end{diagram}
  the top arrow is equal to $X  \alpha  Y$, by $X$ tensor
  (\cel{L}),
  and the left-hand arrow is also equal to $X  \alpha  Y$, by (\cel{R}) tensor $Y$.
  Since $X\alpha Y$ is an isomorphism, it follows that $X \lambda_Y = \rho_X Y$.
\end{dem}
% \begin{dem}
% (This is really a simplicial argument.)  We claim the three
% maps are equal:
% $$
% X  I  \lambda_Y  = X  \alpha  Y = \rho_X
%  I  Y
% $$
% Indeed, the first equality is $X$ tensor (L), and
% the second equality is (R) tensor $Y$.  Now in the commutative 
% square
% \begin{diagram}[w=6ex,h=4.5ex,tight]
% XIIY & \rTo^{\rho_X I Y}  & XIY  \\
% \dTo<{XI\lambda_Y}  &    & \dTo>{X\lambda_Y}  \\
% XIY  & \rTo_{\rho_X Y}  & XY
% \end{diagram}
% we have just shown the first two arrows coincide, and since they are
% isomorphisms, we can cancel them away to
% arrive at the result.
% \end{dem}

\begin{lemma}\label{ass+eq}
  For a unit object $(I,\alpha)$ we have:
  (i) The map $\alpha:II\to I$ is associative.
  (ii) The two functors $X\mapsto IX$ and $X\mapsto XI$ are equivalences.
\end{lemma}

\begin{dem}
  Since $\alpha$ is invertible, associativity amounts to the equation $I \alpha
  = \alpha I$, which follows from the previous proof by setting $X=Y=I$ and
  applying $\cel L$ and $\cel R$ once again.  To see that $\mathbb{L}$ is an
  equivalence, just note that it is isomorphic to the identity via $\lambda$.
\end{dem}

\begin{blanko}{Uniqueness of units.}
  Just as in a semi-monoid a unit element is unique if it exists, one can
  show~\cite{Kock:0507349} that in a semi-monoidal category, any two units are
  uniquely isomorphic.  This statement does not involve $\lambda$ and $\rho$,
  but the proof does: the canonical isomorphism $I \isopil J$ is the composite
  $I \stackrel{\rho_I\inv}{\rTo} IJ \stackrel{\lambda_J}{\rTo} J$.
\end{blanko}

\section{Units in monoidal $2$-categories: definition and main results}
%%%%%%%%%%%%%%%%%%%%%%%%%%%%%%%%%%%%%%%%%%%%%%%%%%
\label{sec:dim2-main}

In this section we set up the necessary terminology and notation, give the main
definition, and state the main results.

\begin{blanko}{$2$-categories.}
  We work in a strict $2$-category $\CC$.  We use the symbol $\#$ to denote
  composition of arrows and horizontal composition of $2$-cells in $\CC$, always
  written from the left to the right, and occasionally decorating the symbol $\#$ by
  the name of the object where the two arrows or $2$-cells are composed.  By an
  {\em equi-arrow} in $\CC$ we understand an arrow $f$ admitting an
  (unspecified) pseudo-inverse, i.e.~an arrow $g$ in the opposite direction such
  that $f \# g$ and $g \# f$ are isomorphic to the respective identity arrows,
  and such that the comparison $2$-cells satisfy the usual triangle equations
  for adjunctions.
  
  We shall make extensive use of arguments with pasting diagrams 
  \cite{Kelly-Street:2cat}.
  Our drawings of $2$-cells should be read from bottom to top, so that for 
  example
  \begin{diagram}[w=3ex,h=3.5ex,tight]
  X    && \rTo^h    && Z    \\
  &\rdTo_f    & \raisebox{4pt}{\cel U}     & \ruTo_g  &  \\
  &    & Y    & &
  \end{diagram}
  denotes $\cel U : f \,\underset Y \#\, g \Rightarrow h$.  The symbol $\comm$ 
  will denote identity $2$-cells.
  
  The few $2$-functors we need all happen to be strict.
  By {\em natural transformation} we always mean pseudo-natural transformation.
  Hence a natural transformation $u: F \Rightarrow G$ between two $2$-functors
  from $\DD$ to $\CC$ is given by an arrow $u_X: FX\to GX$ for each object 
  $X\in \DD$, and an invertible $2$-cell
  \begin{diagram}[w=6ex,h=4.5ex,tight]
  FX & \rTo^{u_X}  & GX  \\
  \dTo<{F(x)}  &  u_x  & \dTo>{G(x)}  \\
  FX'  & \rTo_{u_{X'}}  & GX'
  \end{diagram}
  for each arrow $x: X\to X'$ in $\DD$, subject to the usual compatibility 
  conditions~\cite{Kelly-Street:2cat}.  The modifications we shall need will 
  happen to be invertible.
\end{blanko}

\begin{blanko}{Semi-monoidal $2$-categories.}\label{semimonoidal}
  By {\em semi-monoidal $2$-category} we mean a $2$-category $\CC$ equipped with
  a tensor product, i.e.~an associative $2$-functor
\begin{eqnarray*}
  \tensor : \CC\times\CC & \longrightarrow & \CC  \\
  (X,Y) & \longmapsto & XY ,
\end{eqnarray*}
denoted by plain juxtaposition.  We already assumed $\CC$ to be a strict 
$2$-category, and we also require $\tensor$ to be a strict $2$-functor and to
be strictly associative: $(XY)Z=X(YZ)$.  This is mainly for convenience, to 
keep the focus on unit issues.

Note that the tensor product of two equi-arrows is again an equi-arrow, since
its pseudo-inverse can be taken to be the tensor product of the pseudo-inverses.
\end{blanko}

\begin{blanko}{Semi-monoids.}
  A {\em semi-monoid} in $\CC$ is a triple
  $(X,\alpha,\Aa)$ consisting of an
  object $X$, a multiplication map $\alpha: X X \to X$,
  and an invertible $2$-cell $\Aa$ called the {\em associator},
  $$
  \begin{diagram}[w=6ex,h=4.5ex,tight]
  XXX & \rTo^{\alpha X}  & X X  \\
  \dTo<{X\alpha} &  {\scriptstyle \Aa}  & \dTo>\alpha \\
  X X  & \rTo_{\alpha} & X
  \end{diagram}
$$
  required to satisfy the `pentagon equation':
%%  generate cube \
%%     backobjs XXXX XXX XXX XX \
%%     frontobjs XXX XX XX X \
%%     backarrows {\alpha XX} {XX\alpha} {X\alpha} {\alpha X} \
%%     diagonalarrows  {X\alpha X} {\alpha X} {X\alpha} {\alpha} \
%%     frontarrows {\alpha X} {X\alpha} {\alpha} {\alpha} \
%%     twocells {\comm} {X\Aa} {\Aa X} {\Aa} {\Aa} {\Aa} \
%-------------------------------------output written 2004/11/01 23:02:06 EST----
$$
\begin{diagram}[w=40pt,h=30pt,tight,scriptlabels,objectstyle=\scriptstyle,hug]
  XXXX&\rTo^{\alpha XX}&XXX&&&&\\
  \dTo<{XX\alpha}&\rdTo>{X\alpha X}&\psr{\Aa X}&\rdTo>{\alpha X}&&&\\
  XXX&\psr{X\Aa}&XXX&\rTo^{\alpha X}&XX\\
  &\rdTo_{X\alpha}&\dTo<{X\alpha}&\Aa&\dTo>{\alpha}\\
  &&XX&\rTo_{\alpha}&X
\end{diagram}
\qquad = \qquad
\begin{diagram}[w=40pt,h=30pt,tight,scriptlabels,objectstyle=\scriptstyle,hug]
  XXXX&\rTo^{\alpha XX}&XXX&&&&\\
  \dTo<{XX\alpha}&\comm&\dTo>{X\alpha}&\rdTo>{\alpha X}&&&\\
  XXX&\rTo_{\alpha X}&XX&\psr{\Aa}&XX\\
  &\rdTo_{X\alpha}&\psr{\Aa}&\rdTo<{\alpha}&\dTo>{\alpha}\\
  &&XX&\rTo_{\alpha}&X
\end{diagram}
$$
%-------------------------------------end of output------checksum:0x93C43701----
In the applications, $\alpha$ will be an equi-arrow, and hence we will have
$$
\Aa = \cel A  \underset{XX}{\#}  \alpha
$$
for a some unique invertible
$$
\cel A: X\alpha \Rightarrow \alpha X,
$$ 
which it will more convenient to work with.
% \footnote{The
% Scandinavian letter \AA\ was formed in the middle age, as shorthand 
% for Aa.%
% % the result of a small A creeping up on the shoulder of a bigger A.
% }
In this case, the pentagon equation is equivalent to the more compact equation
%%  generate horizontalCylinder \
%%     startobjects XXXX XXX \
%%     endobjects XXX XX \
%%     straightarrows {\alpha XX} {\alpha X} \
%%     curvedarrows {XX\alpha} {X\alpha X} {X\alpha} {\alpha X} \
%%     twocells {X\cel{A}} {\cel{A}} {\Aa X} {\comm} \
%%     label {shortApentagon}
%-------------------------------------output written 2009/07/12 22:34:20 CEST----
\begin{equation}\label{shortApentagon}
\begin{diagram}[w=6ex,h=6ex,tight]
  XXXX&\rTo^{\alpha XX}&XXX\\
  \leftglob{XX\alpha}\lift{0}{X\cel{A}}\rightglob{X\alpha X}&\shuft{20}{\Aa X}&\rightglob{\alpha X}\\
  XXX&\rTo_{\alpha X}&XX
\end{diagram}
\qquad\qquad = \qquad\qquad
\begin{diagram}[w=6ex,h=6ex,tight]
  XXXX&\rTo^{\alpha XX}&XXX\\
  \leftglob{XX\alpha}&\shuft{-20}{\comm}&\leftglob{X\alpha}\lift{0}{\cel{A}}\rightglob{\alpha X}\\
  XXX&\rTo_{\alpha X}&XX
\end{diagram}
\end{equation}
%-------------------------------------end of output------checksum:0xAD5A6054----
which we shall also make use of.
\end{blanko}

\begin{blanko}{Semi-monoid maps.}\label{semimonoid-map}
  A {\em semi-monoid map} $f: (X,\alpha,\Aa) \to (Y,\beta,\cel{B})$ is
  the data of an arrow $f:X\to Y$ in $\CC$ together with an invertible
  $2$-cell
  \begin{diagram}[w=6ex,h=4.5ex,tight]
  X  X & \rTo^{f f}  & Y Y  \\
  \dTo<\alpha  &  \cel{F}  & \dTo>\beta  \\
  X  & \rTo_f  & Y
  \end{diagram}
  such that this cube commutes:
$$
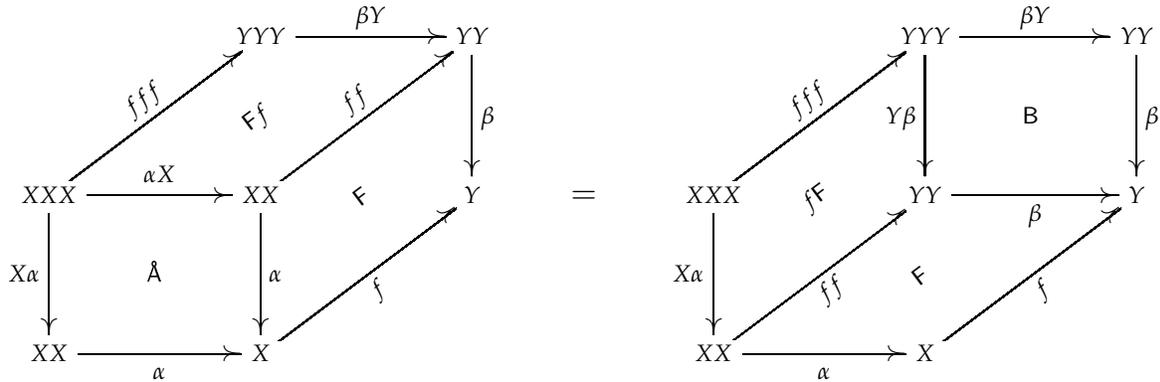
\def\sdf#1{\rotstart{-36.87 rotate}\hbox to0pt {\vsize 0pt \hss\(\scriptstyle #1\)\hss}\rotfinish}
\begin{diagram}[w=40pt,h=30pt,tight,scriptlabels,objectstyle=\scriptstyle,hug]
  &&YYY&\rTo^{\beta Y}& YY \\
  & \ruTo^{fff} &\sdf{\cel{F}f}& \ruTo^{ff} & \dTo>\beta \\
  XXX&\rTo^{\alpha X}&XX&\sdf{\cel F}&Y\\
  \dTo<{X\alpha}&\Aa&\dTo>{\alpha}&\ruTo_f&\\
  XX&\rTo_{\alpha}&X&&
\end{diagram}
\qquad = \qquad
\begin{diagram}[w=40pt,h=30pt,tight,scriptlabels,objectstyle=\scriptstyle,hug]
  &&YYY&\rTo^{\beta Y}&YY\\
  & \ruTo^{fff}&\dTo<{Y \beta}&\cel{B}&\dTo>{\beta}\\
  XXX&\sdf{f\cel F}&YY&\rTo_{\beta}&Y \\
  \dTo<{X\alpha}&\ruTo_{ff}&\sdf{\cel F}&\ruTo_f&\\
  XX&\rTo_{\alpha}&X&&
\end{diagram}
$$
% %%  generate cube \
% %%     backobjs XXX XX XX X \
% %%     frontobjs YYY YY YY Y \
% %%     backarrows {\alpha X} {X \alpha} {\alpha} {\alpha} \
% %%     diagonalarrows  fff ff ff f \
% %%     frontarrows {\beta Y} {Y \beta} {\beta} {\beta} \
% %%     twocells {\Aa} {f\cel{F}} {\cel{F}f} {\cel{F}} {\cel{F}} {\cel{B}} \
% %%     label {}
% %-------------------------------------output written 2004/11/01 22:38:50 EST----
% $$
% \begin{diagram}[w=40pt,h=30pt,tight,scriptlabels,objectstyle=\scriptstyle,hug]
%   XXX&\rTo^{\alpha X}&XX&&&&\\
%   \dTo<{X \alpha}&\rdTo>{fff}&\psr{\cel{F}f}&\rdTo>{ff}&&&\\
%   XX&\psr{f\cel{F}}&YYY&\rTo^{\beta Y}&YY\\
%   &\rdTo_{ff}&\dTo<{Y \beta}&\cel{B}&\dTo>{\beta}\\
%   &&YY&\rTo_{\beta}&Y
% \end{diagram}
% \qquad = \qquad
% \begin{diagram}[w=40pt,h=30pt,tight,scriptlabels,objectstyle=\scriptstyle,hug]
%   XXX&\rTo^{\alpha X}&XX&&&&\\
%   \dTo<{X \alpha}&\Aa&\dTo>{\alpha}&\rdTo>{ff}&&&\\
%   XX&\rTo_{\alpha}&X&\psr{\cel{F}}&YY\\
%   &\rdTo_{ff}&\psr{\cel{F}}&\rdTo<{f}&\dTo>{\beta}\\
%   &&YY&\rTo_{\beta}&Y
% \end{diagram}
% $$
% %-------------------------------------end of output------checksum:0xE78CB7DC----
% (This notion, with $\cel F$ required invertible, is what more precisely could be
% called a {\em strong} semi-monoid map.  It is the only notion we need in this work.)

When $\beta$ is an equi-arrow, the cube equation is equivalent to the simpler 
equation:
% Note that when deriving this simpler equation, it is necessary to
% cancel away to $\cel F$-cells from the right-hand side of the cube
% equation.  This is possible since $\cel F$ is invertible.  Note that
% the two occurrences of the $\cel F$-cell are actually each other's
% inverses.
% 
%%  generate horizontalCylinder \
%%     startobjects XXX XX \
%%     endobjects YYY YY \
%%     straightarrows {fff} {ff} \
%%     curvedarrows {X\alpha} {\alpha X} {Y\beta} {\beta Y} \
%%     twocells {\cel{A}} {\cel{B}} {\cel{F}f} {f\cel{F}} \
%%     label {shortSemimonoidMap}
%-------------------------------------output written 2009/07/12 22:35:33 CEST----
\begin{equation}\label{shortSemimonoidMap}
\begin{diagram}[w=6ex,h=6ex,tight]
  XXX&\rTo^{fff}&YYY\\
  \leftglob{X\alpha}\lift{0}{\cel{A}}\rightglob{\alpha X}&\shuft{20}{\cel{F}f}&\rightglob{\beta Y}\\
  XX&\rTo_{ff}&YY
\end{diagram}
\qquad\qquad = \qquad\qquad
\begin{diagram}[w=6ex,h=6ex,tight]
  XXX&\rTo^{fff}&YYY\\
  \leftglob{X\alpha}&\shuft{-20}{f\cel{F}}&\leftglob{Y\beta}\lift{0}{\cel{B}}\rightglob{\beta Y}\\
  XX&\rTo_{ff}&YY
\end{diagram}
\end{equation}
%-------------------------------------end of output------checksum:0x2C70F638----
which will be useful.
\end{blanko}

\begin{blanko}{Semi-monoid transformations.}\label{semimonoid-transf}
  A {\em semi-monoid transformation} between two parallel semi-monoid maps $(f,\cel{F})$ and
  $(g,\cel{G})$
  is a $2$-cell $\cel{T}: f \Rightarrow g$ in $\CC$ such that this cylinder 
  commutes:
%%  generate verticalCylinder \
%%     startobjects XX YY \
%%     endobjects X Y \
%%     straightarrows {\alpha} {\beta} \
%%     curvedarrows {gg} {ff} {g} {f} \
%%     twocells {\cel{TT}} {\cel{T}} {\cel{F}} {\cel{G}}  \
%%     label {}
%-------------------------------------output written 2009/07/07 14:11:27 CEST----
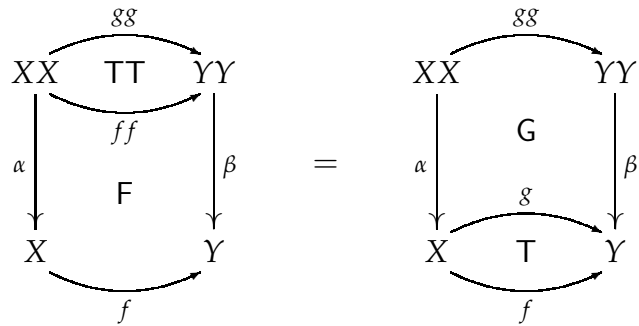
\vspace{24pt}
$$
\begin{diagram}[w=6ex,h=6ex,tight]
  XX&\topglob{gg}\lift{0}{\cel{TT}}\dropglob{ff}&YY\\
  \dTo<{\alpha}&\lift{-12}{\cel{F}}&\dTo>{\beta}\\
  X&\dropglob{f}&Y
\end{diagram}
\qquad = \qquad
\begin{diagram}[w=6ex,h=6ex,tight]
  XX&\topglob{gg}&YY\\
  \dTo<{\alpha}&\lift{12}{\cel{G}}&\dTo>{\beta}\\
  X&\topglob{g}\lift{0}{\cel{T}}\dropglob{f}&Y
\end{diagram}
\vspace{24pt}
$$
%-------------------------------------end of output------checksum:0xD8C17631----
\end{blanko}

% Clearly semi-monoids, semi-monoid maps, and semi-monoid transformations
% constitute a $2$-category, the {\em $2$-category of semi-monoids}.

\begin{lemma}\label{inv}
  Let $f:X \to Y$ be a semi-monoid map.  If $f$ is an equi-arrow (as an arrow in
  $\CC$) with quasi-inverse $g:Y\to X$, then there is a canonical $2$-cell
  $\cel{G}$ such that $(g,\cel{G})$ is a semi-monoid map.
\end{lemma}
\begin{dem}
  The $2$-cell $\cel G$ is defined as the mate~\cite{Kelly-Street:2cat} of
  the $2$-cell $\cel F\inv$.  It is routine to check the cube equation in
  \ref{semimonoid-map}.
\end{dem}
% 
% 
% % 
% The $2$-cell $\cel{G}$ is defined by this pasting diagram:
% \begin{diagram}[w=6ex,h=3ex,tight]
% YY & \rTo^{gg}  & XX  \\
%   &  \luTo_{ff}  & \dLig  \\
%   \dTo<{\beta}&   & XX \\
%   &\cel{F}&\\
%   Y &&\dTo>{\alpha}\\
%   \dLig & \luTo^{f}  & \\
%   Y & \rTo_g & X
% \end{diagram}
% where the two triangles are the adjunction triangles, and $\cel{F}$ is 
% the $2$-cell part of the semi-monoid map $f$.  Checking that that 
% $\cel{G}$ satisfies the cube equation of \ref{semimonoid-map} is routine
% (a bit tedious).
% WE DON'T HAVE TO BECAUSE WE KNOW THAT $I$ AND $J$ ARE UNITS!

\begin{blanko}{Pseudo-idempotents.}\label{idempotent}
  A {\em pseudo-idempotent} is a pair $(I,\alpha)$ where $\alpha:II\to I$ is an
  equi-arrow.  A {\em morphism of pseudo-idempotents} from $(I,\alpha)$ to
  $(J,\beta)$ is a pair $(u,\cel{U})$ consisting of an arrow $u:I\to J$ in $\CC$
  and an invertible $2$-cell
  \begin{diagram}
  II & \rTo^{uu}  & JJ  \\
  \dTo<\alpha  &  \cel{U}  & \dTo>\beta  \\
  I  & \rTo_u  & J .
  \end{diagram}
  If $(u,\cel{U})$ and $(v,\cel{V})$ are morphisms of pseudo-idempotents from
  $(I,\alpha)$ to $(J,\beta)$, a {\em $2$-morphism of pseudo-idempotents} from
  $(u,\cel{U})$ to $(v,\cel{V})$ is a $2$-cell $\cel T: u \Rightarrow v$
  satisfying the cylinder equation of \ref{semimonoid-transf}.  
%   This defines {\em the $2$-category of pseudo-idempotents in $\CC$}.
\end{blanko}

\begin{blanko}{Cancellable objects.}\label{cancellable}
  An object $I$ in $\CC$ is called {\em cancellable} if the two
  $2$-functors $\CC\to\CC$
  \begin{eqnarray*}
    X & \longmapsto & IX\\ 
    X & \longmapsto & XI 
  \end{eqnarray*}
  are fully faithful.  (Fully faithful means that the induced functors on hom
  categories are equivalences.)  A {\em cancellable morphism} between
  cancellable objects $I$ and $J$ is an equi-arrow $u:I \to J$.  (Equivalently
  it is an arrow such that the functors on hom cats defined by tensoring with
  $u$ on either side are equivalences, cf.~\ref{unitmap}.)  A {\em cancellable
  $2$-morphism} between cancellable arrows is any invertible $2$-cell.
%   This defines {\em the category of cancellable objects in $\CC$}.
\end{blanko}

We are now ready for the main definition and the main results.

\begin{blanko}{Units.}\label{main-def}
  A {\em unit object} is by definition a cancellable pseudo-idempotent.
  Hence it is a pair $(I,\alpha)$ consisting
  of an object $I$ and an equi-arrow $\alpha : II \to I$, with the
  property that tensoring with $I$ from either side define fully
  faithful $2$-functors $\CC\to\CC$. 
  
  A {\em morphism} of units is a
  cancellable morphism of
  pseudo-idempotents. In other words, a unit morphism from $(I,\alpha)$ to 
  $(J,\beta)$ is a pair $(u,\cel{U})$ where $u: I \to J$ is an equi-arrow
  and $\cel U$ is an invertible $2$-cell 
  \begin{diagram}
  II & \rTo^{uu}  & JJ  \\
  \dTo<\alpha  &  \cel{U}  & \dTo>\beta  \\
  I  & \rTo_u  & J .
  \end{diagram}

  A {\em $2$-morphism} of units is a cancellable $2$-morphism of
  pseudo-idempotents.  Hence a $2$-morphism from 
  $(u,\cel{U})$ to $(v,\cel{V})$ is a $2$-cell $\cel T: u \Rightarrow v$
  such that
%%  generate verticalCylinder \
%%     startobjects II JJ \
%%     endobjects I J \
%%     straightarrows {\alpha} {\beta} \
%%     curvedarrows {vv} {uu} {v} {u} \
%%     twocells {\cel{TT}} {\cel{T}} {\cel{U}} {\cel{V}}  \
%%     label {}
%-------------------------------------output written 2009/07/15 13:29:36 CEST----
\vspace*{24pt}
$$
\begin{diagram}[w=6ex,h=6ex,tight]
  II&\topglob{vv}\lift{0}{\cel{TT}}\dropglob{uu}&JJ\\
  \dTo<{\alpha}&\lift{-12}{\cel{U}}&\dTo>{\beta}\\
  I&\dropglob{u}&J
\end{diagram}
\qquad = \qquad
\begin{diagram}[w=6ex,h=6ex,tight]
  II&\topglob{vv}&JJ\\
  \dTo<{\alpha}&\lift{12}{\cel{V}}&\dTo>{\beta}\\
  I&\topglob{v}\lift{0}{\cel{T}}\dropglob{u}&J
\end{diagram}
\vspace{24pt}
$$
%-------------------------------------end of output------checksum:0x7866D583----

  This defines the {\em $2$-category of units}.
\end{blanko}
% 
% This notion of unit is very economical: it depends only on the notion of
% equivalence.  
In the next section we'll see how the notion of unit object induces left and
right constraints familiar from standard notions of monoidal $2$-category.  It
will then turn out (Lemmas~\ref{unitmap} and \ref{unit2map}) that unit morphisms
and $2$-morphisms can be characterised as those morphisms and $2$-morphisms
compatible with the left and right constraints.

% It turns out that in that case the
% $2$-functor is actually automatically an equi-$2$-functor.

\begin{blanko*}{Theorem~\ref{thmA} (Associativity).}
  {\em
  Given a unit object $(I,\alpha)$, there is a canonical invertible
  $2$-cell
  \begin{diagram}[w=6ex,h=4.5ex,tight]
  III & \rTo^{\alpha I}  & II  \\
  \dTo<{I\alpha}  &  \Aa  & \dTo>\alpha  \\
  II  & \rTo_\alpha  & I
  \end{diagram}
%   $$
%   I\alpha \stackrel{\cel A}{\Leftrightarrow} \alpha I,
%   $$
  which satisfies the pentagon equation}
%%  generate cube \
%%     backobjs IIII III III II \
%%     frontobjs III II II I \
%%     backarrows {\alpha II} {II\alpha} {I\alpha} {\alpha I} \
%%     diagonalarrows  {I\alpha I} {\alpha I} {I\alpha} {\alpha} \
%%     frontarrows {\alpha I} {I\alpha} {\alpha} {\alpha} \
%%     twocells {\comm} {I\Aa} {\Aa I} {\Aa} {\Aa} {\Aa} \
%%     label pentaA
%-------------------------------------output written 2009/07/10 20:44:00 CEST----
\begin{equation}\label{pentaA}
\begin{diagram}[w=40pt,h=30pt,tight,scriptlabels,objectstyle=\scriptstyle,hug]
  IIII&\rTo^{\alpha II}&III&&&&\\
  \dTo<{II\alpha}&\rdTo>{I\alpha I}&\psr{\Aa I}&\rdTo>{\alpha I}&&&\\
  III&\psr{I\Aa}&III&\rTo^{\alpha I}&II\\
  &\rdTo_{I\alpha}&\dTo<{I\alpha}&\Aa&\dTo>{\alpha}\\
  &&II&\rTo_{\alpha}&I
\end{diagram}
\qquad = \qquad
\begin{diagram}[w=40pt,h=30pt,tight,scriptlabels,objectstyle=\scriptstyle,hug]
  IIII&\rTo^{\alpha II}&III&&&&\\
  \dTo<{II\alpha}&\comm&\dTo>{I\alpha}&\rdTo>{\alpha I}&&&\\
  III&\rTo_{\alpha I}&II&\psr{\Aa}&II\\
  &\rdTo_{I\alpha}&\psr{\Aa}&\rdTo<{\alpha}&\dTo>{\alpha}\\
  &&II&\rTo_{\alpha}&I
\end{diagram}
\end{equation}
%-------------------------------------end of output------checksum:0x4D29F237----
\end{blanko*}
In other words, a unit object is automatically a semi-monoid.
The $2$-cell $\cel A$ is characterised uniquely in \ref{uniqueA}.

\begin{blanko*}{Theorem~\ref{thmB}.}
  {\em
  A unit morphism $(u,\cel{U}): (I,\alpha) \to (J,\beta)$ 
  is automatically a semi-monoid map, when $I$ and $J$ are considered
  semi-monoids in virtue of Theorem~A.}
\end{blanko*}

\begin{blanko*}{Theorem~\ref{thmC} (Contractibility).}
  {\em
  The $2$-category of units in $\CC$ is contractible, if non-empty.
  }
\end{blanko*}
In other words, between any two units there exists a unit morphism, and
between any two parallel unit morphisms there is a unique unit
$2$-morphism.
Theorem~C shows that units objects are unique up to homotopy, so in 
this sense `being unital' is a 
property not a structure.

\bigskip

The proofs of these three theorems rely on the auxiliary structure of left and
right constraints which we develop in the next section, and which also displays
the relation with the classical notion of monoidal $2$-category: in
Section~\ref{sec:comparison} we show that the cancellable-idempotent notion of
unit is equivalent to the notion extracted from the definition of tricategory of
Gordon, Power, and Street~\cite{Gordon-Power-Street}.  This is our 
Theorem~\ref{thmE}.

%%%%%%%%%%%%%%%%%%%%%%%%%%%%%%%%%%%%%%%%%%%%%%%%%%
\section{Left and right actions, and associativity of the unit
(Theorem~\ref{thmA})}
%%%%%%%%%%%%%%%%%%%%%%%%%%%%%%%%%%%%%%%%%%%%%%%%%%
\label{sec:left-right}

Throughout this section we fix a unit object $(I,\alpha)$.
\begin{lemma}\label{LR}
  For each object $X$ there exists pairs $(\lambda_X, \cel L_X)$ and
  $(\rho_X,  \cel R_X)$,
  \begin{eqnarray*}
    \lambda_X : IX \to X, &  &
    \cel L_X : I\lambda_X \Rightarrow \alpha X  \\
    \rho_X : XI \to X,  &  & 
    \cel R_X : X \alpha  \Rightarrow \rho_X I
  \end{eqnarray*}
  where $\lambda_X$ and $\rho_X$ are equi-arrows, and $\cel L_X$ are
  $\cel R_X$ are invertible
  $2$-cells.  
  
  For every such family, there is a unique way to assemble
  the $\lambda_X$ into a natural transformation 
  (this involves
  defining $2$-cells $\lambda_f$ for every arrow $f$ in $\CC$) in such
  a way that $\cel L$ is a natural modification.  Similarly for the
  $\rho_X$ and $\cel R_X$.
%   
%   There exist families of pairs (one for each object $X$),
%   \begin{eqnarray*}
%     (\lambda_X, \cel L_X) & \quad \text{ where } \quad & \lambda_X :
%     IX \to X \text{ is an equi-arrow, and } \cel L_X : I\lambda_X
%     \Leftrightarrow \alpha X \text{ an invertible $2$-cell} \\
%     (\rho_X, \cel R_X) & \quad \text{ where } \quad & \rho_X : XI \to 
%     X \text{ is an equi-arrow, and } \cel R_X : \rho_X I
%     \Leftrightarrow X \alpha \text{ an invertible $2$-cell}
%   \end{eqnarray*}
%   
\end{lemma}

% \begin{lemma}
%   There are natural families of equi-arrows
%   $$
%   \lambda_X : IX \to X
%   $$
%   together with a natural family of invertible $2$-cells
%   $$
%   \cel{L} : I\lambda_X \Leftrightarrow \alpha X .
%   $$
%   
%   Similarly there are natural families of equi-arrows
%   $$
%   \rho_X : XI \to X
%   $$
%   together with natural invertible $2$-cells
%   $$
%   \cel{R} : \rho_X I \Leftrightarrow X \alpha. 
%   $$
%   
%   We refer to these structures as left constraints and right
%   constraints, respectively.
% \end{lemma}

The $\lambda_X$ is an action of $I$ on each $X$, and the $2$-cell
$\cel L_X$ expresses an associativity constraint on this action.
Using these structures we will construct the associator for $\alpha$,
and show it satisfies the pentagon equation. 
Once that is in place we will see
that the actions $\lambda$ and $\rho$ are coherent too (satisfying 
the appropriate pentagon equations).

\bigskip

We shall treat the left action.  The right action is of course
equivalent to establish.

\begin{blanko}{Construction of the left action.}
  Since tensoring with $I$ is a fully faithful $2$-functor, the functor
  $$
  \Hom(IX,X) \to \Hom(IIX,IX)
  $$
  is an equivalence of categories.  In the second category there is
  the canonical object $\alpha  X$.  Hence there is a pseudo 
  pre-image which we denote $\lambda_X : I  X \to X$, together
  with an invertible $2$-cell
  $\cel{L}_X : I  \lambda_X \Rightarrow \alpha  X$:
  \vspace{15pt}
  \begin{diagram}[w=6ex,h=6ex,tight]
  IIX&\topglob{\alpha X}\lift{-2}{\cel L_X}\dropglob{I\lambda_X}&IX
\end{diagram}
\vspace{15pt}

\noindent
  Since $\alpha$ is an equi-arrow, also $\alpha  X$ is equi,
  and since $\cel{L}_X$ is invertible, we conclude that also $I 
  \lambda_X$ is an equi-arrow.  Finally since the $2$-functor `tensoring with
  $I$' is fully faithful, it reflects equi-arrows, so already $\lambda_X$
  is an equi-arrow.  
%   These arguments explain how $\lambda_X$ is
%   constructed.  
%     The $2$-cell $\lambda_f$ for a given arrow $f:X \to Y$ is
%   completely determined by the requirement of compatibility with $\cel{L}$.
\end{blanko}

\begin{blanko}{Naturality.}\label{naturality2}
  A slight variation in the formulation of the construction gives
  directly a natural
  transformation $\lambda$ and a modification $\cel L$:
%   Seeing that these data constitute in fact a natural transformation
%   and that also $\cel{L}$ is natural in $X$, is just a matter of
%   rephrasing the arguments in terms of $2$-functors: 
  Let $\mathbb{L}:
  \CC\to\CC$ denote the $2$-functor `tensoring with $I$ on the left'.
  Since $\mathbb{L}$ is fully faithful, there is an equivalence
  of categories
  $$
  \Nat(\mathbb{L},\Id_\CC) \to \Nat(\mathbb{L}\circ \mathbb{L}, \mathbb{L}) .
  $$
  Now in the second category we have the canonical natural
  transformation whose $X$-component is $\alpha X$ (and with trivial
  components on arrows).  Hence there is a
  pseudo pre-image natural transformation $\lambda : \mathbb{L} \to
  \id_\CC$, together with a modification $\cel L$ whose $X$-component is
  $\cel{L}_X : I \lambda_X \Rightarrow \alpha X$.  
  
  However, we wish to stress the fact that the construction is
  completely
  object-wise.  This fact is of course due to the presence of the
  isomorphism $\cel L_X$: something isomorphic to a natural
  transformation is again natural.  More precisely, to provide the
  $2$-cell data $\lambda_f$ needed to make $\lambda$ into a natural
  transformation, just pull back the $2$-cell data from the natural
  transformation $\alpha X$.
  In detail, we need invertible $2$-cells
%   It is worth spelling out the details of this argument, even were it
%   only to make explicit what $2$-cell data we are talking about.
%   To construct a natural transformation from the data $\lambda_X$
%   is to provide for each arrow $f:X\to Y$ an invertible $2$-cell 
%   \begin{equation}\label{Lambda}
    \begin{diagram}[w=6ex,h=4.5ex,tight]
    I  X  & \rTo^{\lambda_X}  & X  \\
    \dTo<{I f}  &  \lambda_f  & \dTo>f  \\
    I  Y   & \rTo_{\lambda_Y}  & Y
    \end{diagram}
%   \end{equation}
%   respecting   composition of arrows and identity arrows.
%   In detail, the $2$-cell associated to a composite $X\stackrel{f}{\to}
%   Y \stackrel{g}{\to} Z$ is the pasting of the $2$-cells $\lambda_f$
%   and $\lambda_g$, and the $2$-cell associated to an identity arrow is
%   the identity $2$-cell.  Furthermore, the following compatibility
%   with $2$-cells is required.  Given a $2$-cell $\cel T: f \Rightarrow
%   g$, we have this cylinder equation:
% %%  generate horizontalCylinder \
% %%     startobjects IX IY \
% %%     endobjects X Y \
% %%     straightarrows {\lambda_X} {\lambda_Y} \
% %%     curvedarrows {If} {Ig} {f} {g} \
% %%     twocells {I\cel{T}} {\cel{T}} {\lambda_g} {\Lambda_f} \
% %%     label {}
% %-------------------------------------output written 2009/07/13 14:49:33 CEST----
% $$
% \begin{diagram}[w=6ex,h=6ex,tight]
%   IX&\rTo^{\lambda_X}&X\\
%   \leftglob{If}\lift{0}{I\cel{T}}\rightglob{Ig}&\shuft{20}{\lambda_g}&\rightglob{g}\\
%   IY&\rTo_{\lambda_Y}&Y
% \end{diagram}
% \qquad\qquad = \qquad\qquad
% \begin{diagram}[w=6ex,h=6ex,tight]
%   IX&\rTo^{\lambda_X}&X\\
%   \leftglob{If}&\shuft{-20}{\Lambda_f}&\leftglob{f}\lift{0}{\cel{T}}\rightglob{g}\\
%   IY&\rTo_{\lambda_Y}&Y
% \end{diagram}
% $$
% %-------------------------------------end of output------checksum:0xCE4FF297----
% 
To say that the $\cel L_X$ constitute a modification (from $\lambda$
to the identity) is to have this compatibility for every arrow $f:X\to Y$:
%%  generate verticalCylinder \
%%     startobjects IIX IX \
%%     endobjects IIY IY \
%%     straightarrows {IIf} {If} \
%%     curvedarrows {\alpha X} {I\lambda_X} {\alpha Y} {I\lambda_Y} \
%%     twocells {\cel{L}_X} {\cel{L}_Y} {I\lambda_f} {\comm} \
%%     label {}
%-------------------------------------output written 2009/07/24 22:23:53 CEST----
\vspace{24pt}
$$
\begin{diagram}[w=6ex,h=6ex,tight]
  IIX&\topglob{\alpha X}\lift{0}{\cel{L}_X}\dropglob{I\lambda_X}&IX\\
  \dTo<{IIf}&\lift{-12}{I\lambda_f}&\dTo>{If}\\
  IIY&\dropglob{I\lambda_Y}&IY
\end{diagram}
\qquad = \qquad
\begin{diagram}[w=6ex,h=6ex,tight]
  IIX&\topglob{\alpha X}&IX\\
  \dTo<{IIf}&\lift{12}{\comm}&\dTo>{If}\\
  IIY&\topglob{\alpha Y}\lift{0}{\cel{L}_Y}\dropglob{I\lambda_Y}&IY
\end{diagram}
\vspace{24pt}
$$
%-------------------------------------end of output------checksum:0xCFBAC69----
(Here the commutative cell is actually the $2$-cell part of the natural
transformation $\alpha X$.)
Now the point is that each $2$-cell $\lambda_f$ is uniquely defined by
this compatibility: indeed, since the other three $2$-cells in the
diagram are invertible, there is a unique $2$-cell that can fill the
place of $I\lambda_f$, and since $I$ is cancellable this $2$-cell
comes from a unique $2$-cell $\lambda_f$.
The required compatibilities of $\lambda_f$ with composition, identities, 
and $2$-cells now follows from its construction: $\lambda_f$ is
just the translation via $\cel L$ of the identity $2$-cell $\alpha X$.
% which certainly satisfies the compatibility requirements.
\end{blanko}

\begin{blanko}{Uniqueness of the left constraints.}\label{lambda-lambda}
  There may be many choices for $\lambda_X$, and even for a fixed
  $\lambda_X$, there may be many choices for
  $\cel{L}_X$.  However, between any two pairs $(\lambda_X,\cel{L}_X)$ and
  $(\lambda'_X,\cel{L}'_X)$ there is a unique invertible $2$-cell $\cel{U}\loft_X:
  \lambda_X \Rightarrow \lambda'_X$ such that this compatibility
  holds:
%%  generate triangle \
%%     objects {I\lambda_X} {\alpha X} {I \lambda'_X} \
%%     arrows {I \cel{U}\loft_X} {\cel L_X} {\cel L'_X} \
%%     label {}
%-------------------------------------output written 2009/07/14 10:37:47 CEST----
$$
  \xymatrixrowsep{50pt}
  \xymatrixcolsep{42pt}
  \xymatrix @!=0pt {
  I\lambda_X \ar@{=>}[rr]^{I \cel{U}\loft_X} \ar@{=>}[dr]_{\cel L_X} && I \lambda'_X \ar@{=>}[dl]^{\cel L'_X} \\
  & \alpha X &
  }
$$
%-------------------------------------end of output------checksum:0xF95955D3----
  Indeed, this diagram defines uniquely an invertible $2$-cell $I\lambda_X
  \Rightarrow I \lambda_X'$, and since $I$ is cancellable, this
  $2$-cell comes from a unique $2$-cell $\lambda_X \Rightarrow
  \lambda_X'$ which we then call $\cel{U}\loft_X$.

  There is of course a completely analogous statement for right constraints.
\end{blanko}

\begin{blanko}{Construction of the associator.}
  We define $\cel A : I \alpha \Rightarrow \alpha I$
  as the unique $2$-cell satisfying the equation
  \vspace{24pt}
  \begin{equation}\label{A}
  \begin{diagram}[w=56pt,h=42pt,tight]
  IIII & \thintopglob{I\alpha I}\lift{0}{\cel R_I\inv I}\thindropglob{\rho II}  & III  \\
  \thinleftglob{I\alpha I}\lift{0}{I\cel L_I\inv}\thinrightglob{II\lambda}  &  
  \comm
  & \thinleftglob{I\lambda}\lift{0}{\cel{L}_I}\thinrightglob{\alpha I}  \\
  III  & \thintopglob{\rho I}\lift{0}{\cel{R}_I}\thindropglob{I\alpha}  & II
  \end{diagram}
  \qquad\qquad = \qquad
  \begin{diagram}[w=30pt,h=22.5pt,tight,hug]
  IIII &&&&  \\
  &\rdTo^{I\alpha I} &&&& \\
  && III && \\
  &&&\diagonaltopglob{\alpha I}
  \diagonaldropglob{I\alpha}\NWSElabel{\cel A}   & \\
  &&&& II
  \end{diagram}
  \vspace{24pt}
\end{equation}
This definition is meaningful: since $I\alpha I$ is an equi-arrow,
pre-composing with $I\alpha I$ is a $2$-equivalence, hence gives a bijection on 
the level of $2$-cells, so $\cel A$ is 
determined by the left-hand side of the equation.
Note that $\cel A$ is invertible since all the $2$-cells in the
  construction are.  
  
  The associator  $\Aa$ is defined as $\cel A \text{-followed-by-}\alpha$:
  $$
  \Aa \ \df \ \cel A \ \underset{II}{\#} \ \alpha    ,
  $$
  but it will be more convenient to work with $\cel A$.
\end{blanko}

\begin{prop}\label{independence}
  The definition of $\cel{A}$ does not depend on the choices of left constraint
  $(\lambda,\cel{L})$ and right constraint $(\rho,\cel{R})$.
\end{prop}

\begin{dem}
  Write down the left-hand side of \eqref{A} in terms of different left and
  right constraints.  Express these cells in terms of the original $\cel L_I$
  and $\cel R_I$, using the comparison $2$-cells $\cel{U}\loft_I$ and
  $\cel{U}\rught_I$ of \ref{lambda-lambda}.  Finally observe that these
  comparison cells can be moved across the commutative square to cancel each
  other pairwise.
\end{dem}

\begin{blanko}{Uniqueness of $\cel A$.}\label{uniqueA}
  Equation~\eqref{A} may not appear familiar, but it 
  is equivalent to the following
  `pentagon' equation (after post-whiskering with $\alpha$):
%%  generate cube \
%%     backobjs IIII III III II \
%%     frontobjs III II II I \
%%     backarrows {\rho II} {II\lambda} {I\lambda} {\rho I} \
%%     diagonalarrows  {I\alpha I} {\alpha I} {I\alpha} {\alpha} \
%%     frontarrows {\alpha I} {I\alpha} {\alpha} {\alpha} \
%%     twocells {\comm} {(I\cel L)\#(I\alpha)} {(\cel R I)\# (\alpha I)} {\cel L\#\alpha} {\cel R \#\alpha} {\cel A\#\alpha} \
%%     label {penta-LRK}
%-------------------------------------output written 2009/07/07 16:41:03 CEST----
\begin{equation}\label{penta-LRK}
\begin{diagram}[w=40pt,h=30pt,tight,scriptlabels,objectstyle=\scriptstyle,hug]
  IIII&\rTo^{\rho II}&III&&&&\\
  \dTo<{II\lambda}&\rdTo>{I\alpha I}&\psr{(\cel R I)\# (\alpha I)}&\rdTo>{\alpha I}&&&\\
  III&\psr{(I\cel L)\#(I\alpha)}&III&\rTo^{\alpha I}&II\\
  &\rdTo_{I\alpha}&\dTo<{I\alpha}&\cel A\#\alpha&\dTo>{\alpha}\\
  &&II&\rTo_{\alpha}&I
\end{diagram}
\quad = \quad
\begin{diagram}[w=40pt,h=30pt,tight,scriptlabels,objectstyle=\scriptstyle,hug]
  IIII&\rTo^{\rho II}&III&&&&\\
  \dTo<{II\lambda}&\comm&\dTo>{I\lambda}&\rdTo>{\alpha I}&&&\\
  III&\rTo_{\rho I}&II&\psr{\cel L\#\alpha}&II\\
  &\rdTo_{I\alpha}&\psr{\cel R \#\alpha}&\rdTo<{\alpha}&\dTo>{\alpha}\\
  &&II&\rTo_{\alpha}&I
\end{diagram}
\end{equation}
%-------------------------------------end of output------checksum:0x63B82B9----
From this pentagon equation we shall derive the pentagon equation for $\cel A$,
asserted in Theorem~\ref{thmA}.  To this end we need comparison between $\alpha$, 
$\lambda_I$, and $\rho_I$, which we now establish, in analogy with Axiom~(1)
of monoidal category: the left and right constraints
coincide on the unit object, up to a canonical $2$-cell:
\end{blanko}

\begin{lemma}\label{lambda-rho}
  There are unique invertible $2$-cells
  $$
  \rho_I \stackrel{\cel{E}}{\Rightarrow} 
  \alpha \stackrel{\cel{D}}{\Rightarrow} \lambda_I ,
  $$
  such that
  \vspace{30pt}
\begin{equation}
  \label{SAT}
\begin{diagram}[w=7ex,h=6ex,tight]
  III&
  \bigtopglob{\alpha I}\lift{12}{\cel{L}}\lift{-12}{I\cel{D}}\straighthor{I\lambda}\bigdropglob{I\alpha}
  &II
\end{diagram}
\quad = \quad
\begin{diagram}[w=6ex,h=6ex,tight]
  III&
  \topglob{\alpha I}\lift{0}{\cel A}\dropglob{I\alpha}
  &II
\end{diagram}
\quad = \quad
\begin{diagram}[w=7ex,h=6ex,tight]
  III&
  \bigtopglob{\alpha I}\lift{12}{\cel E I}\lift{-12}{\cel R}\straighthor{\rho I}\bigdropglob{I\alpha}
  &II
\end{diagram}
\vspace{30pt}
\end{equation}
% 
% That is, such that these two triangles of $2$-morphisms commute:
% %%  generate twotriangles \
% %%     objects {I\alpha} {I\lambda_I} {\alpha I} \
% %%     arrows {\cel A} {I \cel{D}} {\cel L} \
% %%     objects {I\alpha} {\rho_I I} {\alpha I}   \
% %%     arrows {\cel A} {\cel R} {\cel E I}  \
% %%     label {}
% %-------------------------------------output written 2009/07/07 14:52:41 CEST----
% $$
%   \xymatrixrowsep{50pt}
%   \xymatrixcolsep{42pt}
%   \xymatrix @!=0pt {
%   I\alpha \ar@{=>}[rr]^{\cel A} \ar@{=>}[dr]_{I \cel{D}} && \alpha I \\
%   & I\lambda_I \ar@{=>}[ur]_{\cel L} &
%   }
%   \qquad\qquad
%    \xymatrixrowsep{50pt}
%   \xymatrixcolsep{42pt}
%   \xymatrix @!=0pt {
%   I\alpha \ar@{=>}[rr]^{\cel A} \ar@{=>}[rd]_{\cel R} && \alpha I \\
%   &  \rho_I I \ar@{=>}[ur]_{\cel E I} &
%   }
% $$
% %-------------------------------------end of output------checksum:0x3AC6D5BD----
\end{lemma}
\begin{dem}
  The left-hand equation defines uniquely a $2$-cell $I \alpha \Rightarrow
  I\lambda_I$, and since $I$ is cancellable, this cell comes from a unique
  $2$-cell $\alpha \Rightarrow \lambda_I$ which we then call $\cel{D}$.  Same 
  argument for $\cel{E}$.
\end{dem}

\begin{satz}[Associativity]\label{thmA}
  Given a unit object $(I,\alpha)$, there is a canonical invertible
  $2$-cell
  \begin{diagram}
  III & \rTo^{\alpha I}  & II  \\
  \dTo<{I\alpha}  &  \Aa  & \dTo>\alpha  \\
  II  & \rTo_\alpha  & I
  \end{diagram}
%   $$
%   I\alpha \stackrel{\cel A}{\Leftrightarrow} \alpha I,
%   $$
  which satisfies the pentagon Equation~\eqref{pentaA}.
\end{satz}

\begin{dem}
  On each side of the cube equation \eqref{penta-LRK}, paste the cell $\cel{E}
  II$ on the top, and the cell $II\cel{D}$ on the left.  On the left-hand side
  of the equation we can use Equations~\eqref{SAT} directly, while on the
  right-hand side we first need to move those cells across the commutative
  square before applying \eqref{SAT}.  The result is precisely the pentagon cube
  for $\Aa = \cel A \# \alpha$.
\end{dem}

\begin{blanko}{Coherence of the actions.}
  We have now established that $(I,\alpha,\Aa)$ is a semi-monoid, and
  may observe that the left and right constraints are coherent actions,
  i.e.~that their `associators' $\cel L$ and $\cel R$ satisfy the appropriate
  pentagon equations.  For the left action this equation is:
%%  generate cube \
%%     backobjs IIIX IIX IIX IX \
%%     frontobjs IIX IX IX X \
%%     backarrows {\alpha IX} {II\lambda} {I\lambda} {\alpha X} \
%%     diagonalarrows  {I\alpha X} {\alpha X} {I\lambda} {\lambda} \
%%     frontarrows {\alpha X} {I \lambda} {\lambda} {\lambda} \
%%     twocells {\comm} {(I\cel L) \# (I\lambda)} {\Aa X} {\cel L \#\lambda} {\cel L \#\lambda} {\cel L \#\lambda}
%-------------------------------------output written 2009/07/10 20:50:23 CEST----
$$
\begin{diagram}[w=40pt,h=30pt,tight,scriptlabels,objectstyle=\scriptstyle,hug]
  IIIX&\rTo^{\alpha IX}&IIX&&&&\\
  \dTo<{II\lambda}&\rdTo>{I\alpha X}&\psr{\Aa X}&\rdTo>{\alpha X}&&&\\
  IIX&\psr{(I\cel L) \# (I\lambda)}&IIX&\rTo^{\alpha X}&IX\\
  &\rdTo_{I\lambda}&\dTo<{I \lambda}&\cel L \#\lambda&\dTo>{\lambda}\\
  &&IX&\rTo_{\lambda}&X
\end{diagram}
\qquad = \qquad
\begin{diagram}[w=40pt,h=30pt,tight,scriptlabels,objectstyle=\scriptstyle,hug]
  IIIX&\rTo^{\alpha IX}&IIX&&&&\\
  \dTo<{II\lambda}&\comm&\dTo>{I\lambda}&\rdTo>{\alpha X}&&&\\
  IIX&\rTo_{\alpha X}&IX&\psr{\cel L \#\lambda}&IX\\
  &\rdTo_{I\lambda}&\psr{\cel L \#\lambda}&\rdTo<{\lambda}&\dTo>{\lambda}\\
  &&IX&\rTo_{\lambda}&X
\end{diagram}
$$
%-------------------------------------end of output------checksum:0x27415319----
Establishing this (and the analogous equation for the right action)
is a routine calculation which we omit since we will not 
actually need the result.
We also note that
  the two actions are compatible---i.e.~constitute a two-sided action.
  Precisely this means that there is a canonical invertible $2$-cell
  \begin{diagram}
  I  X  I & \rTo^{\lambda_X  I}  & X  I  \\
  \dTo<{I  \rho_X}  &  \cel{B}  & \dTo>{\rho_X}  \\
  I  X  & \rTo_{\lambda_X}  & X
  \end{diagram}
  This $2$-cell satisfies two pentagon equations, one for $IIXI$ and one for
  $IXII$.
%   , which we do not care to make explicit.
\end{blanko}

\section{Units in the $2$-category of arrows in $\CC$, and Theorem~\ref{thmB}}
%%%%%%%%%%%%%%%%%%%%%%%%%%%%%%%%%%%%%%%%%%%%%%%%%%

\label{sec:catArr}

In this section we prove Theorem~\ref{thmB}, which asserts that
a morphism of units $(u,\cel{U}) : (I,\alpha) \to (J,\beta)$
is automatically a semi-monoid map (with respect to the canonical
associators $\cel A$ and $\cel B$ of the two units).
We have to establish the cube equation of \ref{semimonoid-map},
or in fact the reduced version \eqref{shortSemimonoidMap}.
The strategy to establish Equation~\eqref{shortSemimonoidMap} is to interpret
everything in the $2$-category of arrows of $\CC$.  The key point is to prove
that a morphism of units is itself a unit in the $2$-category of arrows.  Then
we invoke Theorem~\ref{thmA} to get an associator for this unit, and a pentagon equation,
whose short form \eqref{shortApentagon} will be
the sought equation.

\begin{blanko}{The $2$-category of arrows.}
  The {\em $2$-category of arrows} in $\CC$, denoted $\CCa$, is the $2$-category
described as follows.
The objects of $\CCa$ are the
arrows of $\CC$,
$$
X_0 \rTo^x X_1 .
$$
The arrows from $(X_0,X_1,x)$ to $(Y_0,Y_1,y)$ are triples $(f_0,f_1,F)$
where $f_0:X_0 \to Y_0$ and $f_1:X_1\to Y_1$ are arrows in $\CC$
and $F$ is a $2$-cell
\begin{diagram}
X_0 & \rTo^{f_0}  & Y_0  \\
\dTo<x  &  F  & \dTo>y  \\
X_1  & \rTo_{f_1}  & Y_1
\end{diagram}
% for simplicity we only consider invertible $2$-cells.???
If $(g_0,g_1,G)$ is another arrow from $(X_0,X_1,x)$ to $(Y_0,Y_1,y)$,
a $2$-cell from $(f_0,f_1,F)$ to $(g_0,g_1,G)$ is given by a pair $(m_0,m_1)$
where $m_0: f_0\Rightarrow g_0$ and $m_1: f_1 \Rightarrow g_1$ are
$2$-cells in $\CC$ compatible with $F$ and $G$ in the sense that
this cylinder commutes:
%%  generate verticalCylinder \
%%     startobjects X_0 Y_0 \
%%     endobjects X_1 Y_1 \
%%     straightarrows {x} {y} \
%%     curvedarrows {g_0} {f_0} {g_1} {f_1} \
%%     twocells {m_0} {m_1} {F} {G} \
%%     label {}
%-------------------------------------output written 2009/07/13 14:18:32 CEST----
\vspace{24pt}
$$
\begin{diagram}[w=6ex,h=6ex,tight]
  X_0&\topglob{g_0}\lift{0}{m_0}\dropglob{f_0}&Y_0\\
  \dTo<{x}&\lift{-12}{F}&\dTo>{y}\\
  X_1&\dropglob{f_1}&Y_1
\end{diagram}
\qquad = \qquad
\begin{diagram}[w=6ex,h=6ex,tight]
  X_0&\topglob{g_0}&Y_0\\
  \dTo<{x}&\lift{12}{G}&\dTo>{y}\\
  X_1&\topglob{g_1}\lift{0}{m_1}\dropglob{f_1}&Y_1
\end{diagram}
\vspace{24pt}
$$
%-------------------------------------end of output------checksum:0x70E980B8----
Composition of arrows in $\CCa$ is just pasting of squares.
Vertical composition of $2$-cells is just vertical composition of the 
components (the compatibility is guaranteed by pasting of cylinders along
squares), and horizontal composition of $2$-cells is horizontal
composition of the components (compatibility guaranteed by 
pasting along the 
straight sides of the cylinders).
Note that
$\CCa$ inherits a tensor product from $\CC$: this follows from 
functoriality of the tensor product on $\CC$.
% 
% There are obvious functors $\CC\rTo\CCa$ (sending an object to its
% identity arrow) and $\CCa\topile\CC$ (sending an arrow to its
% source, resp.~target).
\end{blanko}

\begin{lemma}\label{map-canc}
  If $I_0$ and $I_1$ are cancellable objects in $\CC$ and $i: I_0 \to 
  I_1$
  is an equi-arrow, then $i$ is cancellable in $\CCa$.
\end{lemma}

\begin{dem}
  We have to show that for given arrows $x:X_0 \to X_1$ and $y: Y_0\to
  Y_1$, the functor
  $$
  \Hom_{\CCa}(x,y) \to \Hom_{\CCa}(ix,iy)
  $$
  defined by tensoring with $i$ on the left is an equivalence of
  categories (the check for tensoring on the right is analogous).
  
  Let us first show that this functor is essentially surjective.
  Let 
  \begin{diagram}
  I_0 X_0 & \rTo^{s_0}  & I_0 Y_0  \\
  \dTo<{ix}  & S   & \dTo>{iy}  \\
  I_1 X_1  & \rTo_{s_1}  & I_1 Y_1
  \end{diagram}
  be an object in $\Hom_{\CCa}(ix,iy)$.  We need to find a square
%   We want to show that this is isomorphic to something of the form
%   $(I_0 k_0, I_1 k_1, i K)$ for some square
  \begin{diagram}
  X_0 & \rTo^{k_0}  & Y_0  \\
  \dTo<{x}  & K   & \dTo>{y}  \\
  X_1  & \rTo_{k_1}  &Y_1
  \end{diagram}
  and an isomorphism $(m_0,m_1)$ from $(s_0,s_1,S)$ to $(I_0 k_0, I_1 k_1, i 
  K)$,
  i.e.~a cylinder
%%  generate verticalCylinder \
%%     startobjects {I_0 X_0} {I_0 Y_0} \
%%     endobjects {I_1 X_1} {I_1 Y_1} \
%%     straightarrows {ix} {iy} \
%%     curvedarrows {I_0 k_0} {s_0} {I_1 k_1} {s_1}  \
%%     twocells {m_0} {m_1} {S} {i K} \
%%     label {}
%-------------------------------------output written 2009/07/14 21:21:53 CEST----
\vspace{24pt}
$$
\begin{diagram}[w=6ex,h=6ex,tight]
  I_0 X_0&\topglob{I_0 k_0}\lift{0}{m_0}\dropglob{s_0}&I_0 Y_0\\
  \dTo<{ix}&\lift{-12}{S}&\dTo>{iy}\\
  I_1 X_1&\dropglob{s_1}&I_1 Y_1
\end{diagram}
\qquad = \qquad
\begin{diagram}[w=6ex,h=6ex,tight]
  I_0 X_0&\topglob{I_0 k_0}&I_0 Y_0\\
  \dTo<{ix}&\lift{12}{i K}&\dTo>{iy}\\
  I_1 X_1&\topglob{I_1 k_1}\lift{0}{m_1}\dropglob{s_1}&I_1 Y_1
\end{diagram}
\vspace{24pt}
$$
%-------------------------------------end of output------checksum:0xD7DA58C4----
Since $I_0$ is a cancellable object, the arrow $s_0$ is isomorphic to $I_0 k_0$
for some $k_0: X_0 \to Y_0$.  Let the connecting
  invertible $2$-cell be denoted $m_0: s_0 \Rightarrow I_0 k_0$.
  Similarly we find $k_1$ and $m_1: s_1 
  \Rightarrow I_1 k_1$.  Since $m_0$ and $m_1$ are invertible, there is a unique
  $2$-cell
  \begin{diagram}
  I_0 X_0 & \rTo^{I_0 k_0}  & I_0 Y_0  \\
  \dTo<{ix}  & T   & \dTo>{iy}  \\
  I_1 X_1  & \rTo_{I_1 k_1}  & I_1 Y_1
  \end{diagram}
  that can take the place of $iK$ in the cylinder 
  equation; it remains to see that $T$ is of the form $iK$ for some $K$.
  But this follows since the map
  \begin{eqnarray}
    {\rm 2Cell}_{\CC}(k_0\# y, x\# k_1) & \longrightarrow & {\rm 2Cell}_{\CC}(i(k_0\# y), 
    i(x\# k_1))  \notag  \\
    K & \longmapsto & iK \label{iK}
  \end{eqnarray}
  is a bijection.  Indeed, the map factors as `tensoring with $I_0$ on the left'
  followed by `post-composing with $i Y_1$'; the first is a bijection since
  $I_0$ is cancellable, the second is a bijection since $i$ (and hence $i Y_1$) is 
  an equi-arrow).
%   , hence on the level of $2$-cells we have
%   a bijection.  So the $2$-cell $T$ is uniquely given as $iK$ for some
%   $2$-cell $K$.  The cylinder is now the witness that $(m_0,m_1)$ is
%   a $2$-cell in $\CCa$ between our original object and the one in
%   the image.  Hence the functor is essentially surjective.

  Now for the fully faithfulness of 
  $\Hom_{\CCa}(x,y) \to \Hom_{\CCa}(ix,iy)$.
  Fix two objects in the left-hand category, $P$ and $Q$:
  $$
  \begin{diagram}
  X_0 & \rTo^{p_0}  & Y_0  \\
  \dTo<x  & P   & \dTo>y  \\
  X_1  & \rTo_{p_1}  & Y_1
  \end{diagram}
  \qquad
  \begin{diagram}
  X_0 & \rTo^{q_0}  & Y_0  \\
  \dTo<x  & Q   & \dTo>y  \\
  X_1  & \rTo_{q_1}  & Y_1
  \end{diagram}
$$
The arrows from $P$ to $Q$ are pairs $(m_0,m_1)$ consisting of
$$
m_0: p_0 \Rightarrow q_0 \qquad m_1: p_1 \Rightarrow q_1
$$
cylinder-compatible with the $2$-cells $P$ and $Q$.
The image of these two objects are
$$
\begin{diagram}
I_0 X_0 & \rTo^{I_0 p_0}  & I_0 Y_0  \\
\dTo<{ix}  & iP   & \dTo>{iy}  \\
I_1 X_1  & \rTo_{I_1 p_1}  & I_1 Y_1
\end{diagram}
\qquad
\begin{diagram}
I_0 X_0 & \rTo^{I_0 q_0}  & I_0 Y_0  \\
\dTo<{ix}  & iQ   & \dTo>{iy}  \\
I_1 X_1  & \rTo_{I_1 q_1}  & I_1 Y_1
\end{diagram}
$$
The possible $2$-cells from $iP$ to $iQ$ are
pairs $(n_0,n_1)$ consisting of 
$$
n_0: I_0 p_0 \Rightarrow I_0 q_0 \qquad n_1: I_1 p_1 \Rightarrow I_1 q_1
$$
cylinder-compatible with the $2$-cells $iP$ and $iQ$.  Now since $I_0$
is cancellable, every $2$-cell $n_0$ like this is uniquely of the form
$I_0 n_0$ for some $n_0$.  Hence there is a bijection between the
possible $m_0$ and the possible $n_0$.  Similarly for $m_1$ and $n_1$.
So there is a bijection between pairs $(m_0, m_1)$ and pairs $(n_0,
n_1)$.  Now by functoriality of tensoring with $i$, all images of
compatible $(m_0,m_1)$ are again compatible.  It remains to rule out 
the possibility that
some $(n_0,n_1)$ pair could be compatible without $(m_0, m_1)$ being so,
but this follows again from the argument that `tensoring with $i$ on the left'
is a  bijection on hom sets, just like argued for \eqref{iK}.
\end{dem}

\begin{lemma}\label{f0f1-equi}
  An arrow in $\CCa$,
  \begin{diagram}
  X_0 & \rTo^{f_0}  & Y_0  \\
  \dTo<x  &  F  & \dTo>y  \\
  X_1  & \rTo_{f_1}  & Y_1
  \end{diagram}
  is an equi-arrow in $\CCa$ if the components 
  $f_0$ and $f_1$ are equi-arrows in $\CC$ and $F$ is invertible.
\end{lemma}
% In fact I don't know about the `only-if' part.  But we only need the
% `if' part, namely that if the components are equi then so is the whole
% thing.

\begin{dem}
  We can construct an explicit quasi-inverse by choosing quasi-inverses 
  to the components.%
\end{dem}

\begin{cor}
  If $(I_0,\alpha_0)$ and $(I_1, \alpha_1)$ are units in $\CC$, and 
  $(u,\cel{U}): I_0 \to I_1$ is a unit map between them, then
  $$
  u: I_0 \to I_1
  $$
  is a unit object in $\CCa$ with structure map
  \begin{diagram}
  I_0 I_0 & \rTo^{\alpha_0}  & I_0  \\
  \dTo<{uu}  & \cel{U}\inv   & \dTo>u \\
  I_1 I_1  & \rTo_{\alpha_1}  & I_1 .
  \end{diagram}
\end{cor}

\begin{dem}
  The object $u$ is cancellable by Lemma~\ref{map-canc}, and the morphism
  $(\alpha_0,\alpha_1,\cel{U}\inv)$ from $uu$ to $u$ is an equi-arrow by 
  Lemma~\ref{f0f1-equi}.
\end{dem}

\begin{satz}\label{thmB}
  Let $(I_0,\alpha_0)$ and $(I_1,\alpha_1)$ be units, with canonical
  associators $\cel A_0$ and $\cel A_1$, respectively.  If $(u,\cel{U})$ is a unit map
  from $I_0$ to $I_1$ then it is automatically a semi-monoid map.  That is,
%%  generate horizontalCylinder \
%%     startobjects {I_0 I_0 I_0} {I_0 I_0} \
%%     endobjects {I_1 I_1 I_1} {I_1 I_1} \
%%     straightarrows {uuu} {uu} \
%%     curvedarrows {I_0\alpha_0} {\alpha_0 I_0} {I_1\alpha_1} {\alpha_1 I_1} \
%%     twocells {\cel A_0} {\cel A_1} {\cel{U}u} {u\cel{U}} \
%%     label {}
%-------------------------------------output written 2009/07/13 14:56:31 CEST----
$$
\begin{diagram}[w=6ex,h=6ex,tight]
  I_0 I_0 I_0&\rTo^{uuu}&I_1 I_1 I_1\\
  \leftglob{I_0\alpha_0}\lift{0}{\cel A_0}\rightglob{\alpha_0 I_0}&\shuft{20}{\cel{U}u}&\rightglob{\alpha_1 I_1}\\
  I_0 I_0&\rTo_{uu}&I_1 I_1
\end{diagram}
\qquad\qquad = \qquad\qquad
\begin{diagram}[w=6ex,h=6ex,tight]
  I_0 I_0 I_0&\rTo^{uuu}&I_1 I_1 I_1\\
  \leftglob{I_0\alpha_0}&\shuft{-20}{u\cel{U}}&\leftglob{I_1\alpha_1}\lift{0}{\cel A_1}\rightglob{\alpha_1 I_1}\\
  I_0 I_0&\rTo_{uu}&I_1 I_1
\end{diagram}
$$
%-------------------------------------end of output------checksum:0x6D2DFD9B----
\end{satz}

\begin{dem}
  By the previous Corollary, $(u,\cel{U}\inv)$ is a unit object in $\CCa$.
  Hence there is a canonical associator
  $$
  \cel B : u\cel{U}\inv \Leftrightarrow \cel{U}\inv u .
  $$
  By definition of $2$-cells in $\CCa$, this is a pair of $2$-cells in $\CC$
  $$
  \cel B_0 : I_0 \alpha_0 \Rightarrow \alpha_0 I_0
  \qquad
  \cel B_1 : I_1 \alpha_1 \Rightarrow \alpha_1 I_1 ,
  $$
  fitting the cylinder equation
%%  generate verticalCylinder \
%%     startobjects {I_0 I_0 I_0} {I_0 I_0} \
%%     endobjects {I_1 I_1 I_1} {I_1 I_1} \
%%     straightarrows {uuu} {uu} \
%%     curvedarrows {\alpha_0 I_0} {I_0 \alpha_0} {\alpha_1 I_1} {I_1 \alpha_1} \
%%     twocells {\cel B_0} {\cel B_1} {u \cel{U}\inv} {\cel{U}\inv u} \
%%     label {}
%-------------------------------------output written 2009/07/13 14:57:53 CEST----
\vspace{24pt}
$$
\begin{diagram}[w=6ex,h=6ex,tight]
  I_0 I_0 I_0&\topglob{\alpha_0 I_0}\lift{0}{\cel B_0}\dropglob{I_0 \alpha_0}&I_0 I_0\\
  \dTo<{uuu}&\lift{-12}{u \cel{U}\inv}&\dTo>{uu}\\
  I_1 I_1 I_1&\dropglob{I_1 \alpha_1}&I_1 I_1
\end{diagram}
\qquad = \qquad
\begin{diagram}[w=6ex,h=6ex,tight]
  I_0 I_0 I_0&\topglob{\alpha_0 I_0}&I_0 I_0\\
  \dTo<{uuu}&\lift{12}{\cel{U}\inv u}&\dTo>{uu}\\
  I_1 I_1 I_1&\topglob{\alpha_1 I_1}\lift{0}{\cel B_1}\dropglob{I_1 \alpha_1}&I_1 I_1
\end{diagram}
\vspace{24pt}
$$
%-------------------------------------end of output------checksum:0x521021CC----
  This is precisely the cylinder diagram we are looking for---provided we can
  show that $\cel B_0 = \cel A_0$ and $\cel B_1 = \cel A_1$.  But this is a
  consequence of the characterising property of the associator of a unit: first
  note that as a unit object in $\CCa$, $u$ induces left and right constraints: for each
  object $x: X_0 \to X_1$ in $\CCa$ there is a left action of the unit $u$, and 
  this left action will induce a left action of $(I_0,\alpha_0)$ on $X_0$
  and a left action of $(I_1,\alpha_1)$ on $X_1$ (the ends of the cylinders).
  Similarly there is a right action of $u$ which induces right actions at the 
  ends of the cylinder.  Now the unique $\cel B$ that exists as
  associator for the unit object $u$ compatible with the left and right
  constraints induces $\cel B_0$ and $\cel B_1$ at the ends of the cylinder, and
  these will of course be compatible with the induced left and right
  constraints. Hence, by uniqueness of associators compatible with left and
  right constraints, these induced associators $\cel B_0$ and $\cel B_1$
  must coincide with $\cel A_0$ and
  $\cel A_1$.  Note that this does not dependent on choice of left and right 
  constraints, cf.~Proposition~\ref{independence}.
\end{dem}

%%%%%%%%%%%%%%%%%%%%%%%%%%%%%%%%%%%%%%%%%%%%%%%%%%
\section{Contractibility of the space of weak units (Theorem~\ref{thmC})}
%%%%%%%%%%%%%%%%%%%%%%%%%%%%%%%%%%%%%%%%%%%%%%%%%%
\label{sec:contractibility}

The goal of this section is to prove Theorem~\ref{thmC}, which asserts that 
the $2$-category of units in $\CC$ is contractible if non-empty.
First we describe the unit morphisms and unit $2$-morphisms in terms of
compatibility with left and right constraints.  This will show that there 
are not too many $2$-cells. Second we use the left
and right constraints to connect any two units.

The following lemma shows that just as the single arrow $\alpha$ induces all the
$\lambda_X$ and $\rho_X$, 
% and the single $2$-cell $\cel A$ induces all the $\cel
% L_X$ and $\cel R_X$, in the same way, 
the single $2$-cell $\cel{U}$ of a unit
map induces families $\cel{U}\loft_X$ and $\cel{U}\rught_X$ expressing
compatibility with $\lambda_X$ and $\rho_X$.

\begin{lemma}\label{unitmap}
  Let $(I,\alpha)$ and $(J,\beta)$ be units, and let
  $(u,\cel{U})$ be a morphism of pseudo-idempotents from $(I,\alpha)$ to $(J,\beta)$.
  The following are equivalent.
\begin{punkt-i}
  \item $u$ is an equi-arrow (i.e.~$u$ is a morphism of units).
  
  \item $u$ is left cancellable,  i.e.~tensoring with $u$ on the left
   is an equivalence of categories $\Hom(X,Y) \to \Hom(IX,JY)$.

   \item[(ii')] $u$ is right cancellable,  i.e.~tensoring with $u$ on the right
   is an equivalence of categories $\Hom(X,Y) \to \Hom(XI,YJ)$.

  \item For fixed left actions $(\lambda_X,\cel L_X)$ for the unit $(I,\alpha)$ and
  $(\ell_X,\cel L'_X)$ for the unit $(J,\beta)$,
  there is a unique invertible $2$-cell $\cel{U}\loft_X$, natural in $X$:
  \begin{diagram}
  IX & \rTo^{uX}  & JX  \\
  \dTo<{\lambda_X}  &  \cel{U}\loft_X  & \dTo>{\ell_X}  \\
  X  & \rTo_X  & X
  \end{diagram}
such that this compatibility holds:
%%  generate horizontalCylinder \
%%     startobjects IIX IX \
%%     endobjects JJX JX \
%%     straightarrows {uuX} {uX} \
%%     curvedarrows {I\lambda_X} {\alpha X} {J\ell_X} {\beta X} \
%%     twocells {\cel{L}_X} {\cel{L}'_X} {\cel{U}X} {u\cel{U}\loft_X} \
%%     label {P}
%-------------------------------------output written 2009/07/15 14:07:48 CEST----
\begin{equation}\label{P}
\begin{diagram}[w=6ex,h=6ex,tight]
  IIX&\rTo^{uuX}&JJX\\
  \leftglob{I\lambda_X}\lift{0}{\cel{L}_X}\rightglob{\alpha X}&\shuft{20}{\cel{U}X}&\rightglob{\beta X}\\
  IX&\rTo_{uX}&JX
\end{diagram}
\qquad\qquad = \qquad\qquad
\begin{diagram}[w=6ex,h=6ex,tight]
  IIX&\rTo^{uuX}&JJX\\
  \leftglob{I\lambda_X}&\shuft{-20}{u\cel{U}\loft_X}&\leftglob{J\ell_X}\lift{0}{\cel{L}'_X}\rightglob{\beta X}\\
  IX&\rTo_{uX}&JX
\end{diagram}
\end{equation}
%-------------------------------------end of output------checksum:0xD03D041----
% 
% 
% 
% 
%   such that this cube equation holds:
%   
% %%  generate cube \
% %%     backobjs IIX IX JJX JX \
% %%     frontobjs IX X JX X \
% %%     backarrows {\alpha X} {uuX} {uX} {\beta X} \
% %%     diagonalarrows  {I\lambda} {\lambda} {I\ell} {\ell} \
% %%     frontarrows {\lambda} {uX} {X} {\ell} \
% %%     twocells {\cel{U}X} {u\cel{U}\loft} {\cel{L}} {\cel{U}\loft} {\cel{L}} {\cel{U}\loft}
% %-------------------------------------output written 2004/09/17 22:25:56 EDT----
% $$
% \begin{diagram}[w=40pt,h=30pt,tight,scriptlabels,objectstyle=\scriptstyle,hug]
%   IIX&\rTo^{\alpha X}&IX&&&&\\
%   \dTo<{ffX}&\rdTo>{I\lambda}&\psr{\cel{L}}&\rdTo>{\lambda}&&&\\
%   JJX&\psr{f\cel{U}\loft}&IX&\rTo^{\lambda}&X\\
%   &\rdTo_{I\ell}&\dTo<{fX}&\cel{U}\loft&\dTo>{X}\\
%   &&JX&\rTo_{\ell}&X
% \end{diagram}
% \qquad = \qquad
% \begin{diagram}[w=40pt,h=30pt,tight,scriptlabels,objectstyle=\scriptstyle,hug]
%   IIX&\rTo^{\alpha X}&IX&&&&\\
%   \dTo<{ffX}&\cel{F}X&\dTo>{fX}&\rdTo>{\lambda}&&&\\
%   JJX&\rTo_{\beta X}&JX&\psr{\cel{U}\loft}&X\\
%   &\rdTo_{I\ell}&\psr{\cel{L}}&\rdTo<{\ell}&\dTo>{X}\\
%   &&JX&\rTo_{\ell}&X
% \end{diagram}
% $$
% %-------------------------------------end of output------checksum:0x8B85CD29----
% and similarly, there is a unique $2$-cell $\cel{U}\rught$ for fixed right actions.

  \item[(iii')] For fixed right actions $(\rho_X,\cel R_X)$ for the unit $(I,\alpha)$ and
  $(r_X,\cel R'_X)$ for the unit $(J,\beta)$,
  there is a unique invertible $2$-cell $\cel{U}\rught_X$, natural in $X$:
  \begin{diagram}
  XI & \rTo^{Xu}  & XJ  \\
  \dTo<{\rho_X}  &  \cel{U}\rught_X  & \dTo>{r_X}  \\
  X  & \rTo_X  & X
  \end{diagram}
such that this compatibility holds:
%%  generate horizontalCylinder \
%%     startobjects XII XI \
%%     endobjects XJJ XJ \
%%     straightarrows {Xuu} {Xu} \
%%     curvedarrows {X\alpha} {\rho_X I} {X\beta} {r_X J} \
%%     twocells {\cel{R}_X} {\cel{R}'_X} {\cel{U}\rught_X u} {X\cel{U}} \
%%     label {Q}
%-------------------------------------output written 2009/07/15 14:04:54 CEST----
\begin{equation}\label{Q}
\begin{diagram}[w=6ex,h=6ex,tight]
  XII&\rTo^{Xuu}&XJJ\\
  \leftglob{X\alpha}\lift{0}{\cel{R}_X}\rightglob{\rho_X I}&\shuft{20}{\cel{U}\rught_X u}&\rightglob{r_X J}\\
  XI&\rTo_{Xu}&XJ
\end{diagram}
\qquad\qquad = \qquad\qquad
\begin{diagram}[w=6ex,h=6ex,tight]
  XII&\rTo^{Xuu}&XJJ\\
  \leftglob{X\alpha}&\shuft{-20}{X\cel{U}}&\leftglob{X\beta}\lift{0}{\cel{R}'_X}\rightglob{r_X J}\\
  XI&\rTo_{Xu}&XJ
\end{diagram}
\end{equation}
%-------------------------------------end of output------checksum:0x618FBD57----

\end{punkt-i}
\end{lemma}

% \begin{blanko}{The category of units.}
% A {\em morphism of left (resp.~right) units} is a semi-monoid map $f:I
% \to J$ satisfying the conditions of the lemma above.
% 
% A {\em $2$-morphism of left (resp.~right) units} is any semi-monoid transformation.
% This defines the $2$-category $\UU$ of units in $\CC$.
% \end{blanko}

\begin{dem}
  (i) implies (ii):   `tensoring with $u$' can be done in two steps:
  given an arrow $X \to Y$, first tensor with $I$ to get $IX \to IY$,
  and then post-compose with $uY$ to get $IX \to JY$.  The first step
  is an equivalence because $I$ is a unit, and the second step is an
  equivalence because $u$ is an equi-arrow.
  
  (ii) implies (iii): 
%   In fact we will establish the following shortened version of the cube 
%   equation:
%     
% %%  generate verticalCylinder \
% %%     startobjects IIX IX \
% %%     endobjects JJX JX \
% %%     straightarrows {ffX} {fX} \
% %%     curvedarrows {\alpha X} {I\lambda} {\beta X} {J\ell} \
% %%     twocells {\cel{L}} {\cel{L}} {f\cel{U}\loft} {F\cel{X}} \
% %%     label {P}
% %-------------------------------------output written 2004/05/24 18:58:32 EDT----
% \vspace{24pt}
% \begin{equation}\label{P}
% \begin{diagram}[w=6ex,h=6ex,tight]
%   IIX&\topglob{\alpha X}\lift{0}{\cel{L}}\dropglob{I\lambda}&IX\\
%   \dTo<{ffX}&\lift{-12}{f\cel{U}\loft}&\dTo>{fX}\\
%   JJX&\dropglob{J\ell}&JX
% \end{diagram}
% \qquad = \qquad
% \begin{diagram}[w=6ex,h=6ex,tight]
%   IIX&\topglob{\alpha X}&IX\\
%   \dTo<{ffX}&\lift{12}{F\cel{X}}&\dTo>{fX}\\
%   JJX&\topglob{\beta X}\lift{0}{\cel{L}}\dropglob{J\ell}&JX
% \end{diagram}
% \vspace{24pt}
% \end{equation}
% %-------------------------------------end of output------checksum:0x3DAF3AC8----
% (The cube equation follows by pasting with a $\cel{U}\loft$-cell on the 
% right.)
In Equation~(\ref{P}), the $2$-cell labelled $u\cel{U}\loft_X$ is uniquely 
defined by the three other cells, and it is invertible since the three other 
cells are.  Since tensoring with $u$ on the 
left is an equivalence, this cell comes from a unique invertible cell 
$\cel{U}\loft_X$,
justifying the label $u\cel{U}\loft_X$.

(iii) implies (i): The invertible $2$-cell $\cel{U}\loft_X$ shows that $uX$ is
  isomorphic to an equi-arrow, and hence is an equi-arrow itself.  Now take $X$
  to be a right cancellable object (like for example $I$) and conclude that
  already $u$ is an equi-arrow.
  
  Finally, the equivalence 
  (i)\impl(ii')\impl(iii')\impl(i) is completely analogous.
\end{dem}

Note that for $(u,\cel U)$ the identity morphism on $(I,\alpha)$, we recover
Observation~\ref{lambda-lambda}.
  
\begin{lemma}\label{unit2map}
  Let $(I,\alpha)$ and $(J,\beta)$ be units; let $(u,\cel{U})$ and $(v,\cel{V})$
  be morphisms of pseudo-idempotents from $I$ to $J$; and consider a
  $2$-cell $\cel T: u\Rightarrow v$.  Then
  the following are equivalent.
  \begin{punkt-i}
    \item $\cel T$ is an invertible $2$-morphism of pseudo-idempotents.
    
    \item $\cel T$ is a left cancellable $2$-morphism of pseudo-idempotents
    (i.e., induces a bijection on hom sets (of hom cats) by tensoring with 
    $\cel T$ from the left).
    
    \item[(ii')] $\cel T$ is a right cancellable $2$-morphism of pseudo-idempotents
    (i.e., induces a bijection on hom sets (of hom cats) by tensoring with 
    $\cel T$ from the right).
    
    \item For fixed left actions $(\lambda_X, \cel L_X)$ for $(I,\alpha)$
    and $(\ell_X,\cel L'_X)$ for $(J,\beta)$,
    with induced canonical $2$-cells $\cel{U}\loft_X$ and
    $\cel{V}\loft_X$ as in \ref{unitmap}, we have:
%%  generate verticalCylinder \
%%     startobjects IX JX \
%%     endobjects X X \
%%     straightarrows {\lambda_X} {\ell_X} \
%%     curvedarrows {vX} {uX} {X} {X} \
%%     twocells {\cel T X} {\comm} {\cel{U}\loft_X} {\cel{V}\loft_X} \
%%     label {TXP}
%-------------------------------------output written 2009/07/14 11:09:05 CEST----
\vspace{24pt}
\begin{equation}\label{TXP}
\begin{diagram}[w=6ex,h=6ex,tight]
  IX&\topglob{vX}\lift{0}{\cel T X}\dropglob{uX}&JX\\
  \dTo<{\lambda_X}&\lift{-12}{\cel{U}\loft_X}&\dTo>{\ell_X}\\
  X&\dropglob{X}&X
\end{diagram}
\qquad = \qquad
\begin{diagram}[w=6ex,h=6ex,tight]
  IX&\topglob{vX}&JX\\
  \dTo<{\lambda_X}&\lift{12}{\cel{V}\loft_X}&\dTo>{\ell_X}\\
  X&\topglob{X}\lift{0}{\comm}\dropglob{X}&X
\end{diagram}
\vspace{24pt}
\end{equation}
%-------------------------------------end of output------checksum:0xC059827E----

    \item[(iii')] For fixed right actions $(\rho_X,\cel R_X)$ for $(I,\alpha)$
    and $(r_X,\cel R'_X)$ for $(J,\beta)$, with induced canonical $2$-cells $\cel{U}\rught_X$ and
    $\cel{V}\rught_X$ as in \ref{unitmap}, we have:
%%  generate verticalCylinder \
%%     startobjects XI XJ \
%%     endobjects X X \
%%     straightarrows {\rho_X} {r_X} \
%%     curvedarrows {Xv} {Xu} {X} {X} \
%%     twocells {X\cel T} {\comm} {\cel{U}\rught_X} {\cel{V}\rught_X} \
%%     label {TXQ}
%-------------------------------------output written 2009/07/15 14:13:51 CEST----
\vspace{24pt}
\begin{equation}\label{TXQ}
\begin{diagram}[w=6ex,h=6ex,tight]
  XI&\topglob{Xv}\lift{0}{X\cel T}\dropglob{Xu}&XJ\\
  \dTo<{\rho_X}&\lift{-12}{\cel{U}\rught_X}&\dTo>{r_X}\\
  X&\dropglob{X}&X
\end{diagram}
\qquad = \qquad
\begin{diagram}[w=6ex,h=6ex,tight]
  XI&\topglob{Xv}&XJ\\
  \dTo<{\rho_X}&\lift{12}{\cel{V}\rught_X}&\dTo>{r_X}\\
  X&\topglob{X}\lift{0}{\comm}\dropglob{X}&X
\end{diagram}
\vspace{24pt}
\end{equation}
%-------------------------------------end of output------checksum:0x2F025384----
\end{punkt-i}
\end{lemma}

\begin{dem}
  It is obvious that (i) implies (ii).
  Let us prove that (ii) implies (iii), so assume that tensoring with
  $\cel T$ on the left defines a bijection on the level of $2$-cells.
  Start with the cylinder diagram for compatibility of tensor
  $2$-cells (cf.~\ref{semimonoid-transf}).  Tensor this diagram with $X$ on
  the right to get
%%  generate verticalCylinder \
%%     startobjects IIX JJX \
%%     endobjects IX JX \
%%     straightarrows {\alpha X} {\beta X} \
%%     curvedarrows {vvX} {uuX} {vX} {uX}\
%%     twocells {\cel{TT}X} {\cel{T}X} {\cel{U}X} {\cel{V}X} \
%-------------------------------------output written 2009/07/13 14:41:35 CEST----
\vspace{24pt}
$$
\begin{diagram}[w=6ex,h=6ex,tight]
  IIX&\topglob{vvX}\lift{0}{\cel{TT}X}\dropglob{uuX}&JJX\\
  \dTo<{\alpha X}&\lift{-12}{\cel{U}X}&\dTo>{\beta X}\\
  IX&\dropglob{uX}&JX
\end{diagram}
\qquad = \qquad
\begin{diagram}[w=6ex,h=6ex,tight]
  IIX&\topglob{vvX}&JJX\\
  \dTo<{\alpha X}&\lift{12}{\cel{V}X}&\dTo>{\beta X}\\
  IX&\topglob{vX}\lift{0}{\cel{T}X}\dropglob{uX}&JX
\end{diagram}
\vspace{24pt}
$$
%-------------------------------------end of output------checksum:0x679423AB----
% where $\cel{U}$ and $\cel{V}$ are the $2$-cell part of $f$ and $g$.
On each side of this equation, paste an $\cel L_X$ along $\alpha X$,
apply Equation~\eqref{P} on each side, and cancel the $\cel L'_X$
that appear on the other side of the square.
The resulting diagram
%%  generate verticalCylinder \
%%     startobjects IIX JJX \
%%     endobjects IX JX \
%%     straightarrows {I\lambda_X} {J \ell_X} \
%%     curvedarrows {vvX} {uuX} {vX} {uX}\
%%     twocells {\cel{TT}X} {\cel{T}X} {u \cel{U}\loft_X} {v\cel{V}\loft_X} \
%-------------------------------------output written 2009/07/14 11:10:19 CEST----
\vspace*{24pt}
$$
\begin{diagram}[w=6ex,h=6ex,tight]
  IIX&\topglob{vvX}\lift{0}{\cel{TT}X}\dropglob{uuX}&JJX\\
  \dTo<{I\lambda_X}&\lift{-12}{u \cel{U}\loft_X}&\dTo>{J \ell_X}\\
  IX&\dropglob{uX}&JX
\end{diagram}
\qquad = \qquad
\begin{diagram}[w=6ex,h=6ex,tight]
  IIX&\topglob{vvX}&JJX\\
  \dTo<{I\lambda_X}&\lift{12}{v\cel{V}\loft_X}&\dTo>{J \ell_X}\\
  IX&\topglob{vX}\lift{0}{\cel{T}X}\dropglob{uX}&JX
\end{diagram}
\vspace{24pt}
$$
%-------------------------------------end of output------checksum:0x2B1B9F10----
is the tensor product of $\cel T$
% %%  generate verticalCylinder \
% %%     startobjects I J \
% %%     endobjects I J \
% %%     straightarrows {I} {J} \
% %%     curvedarrows {v} {u} {v} {u} \
% %%     twocells {\cel{T}} {\cel{T}} {\comm} {\comm} \
% %%     label {}
% %-------------------------------------output written 2009/07/09 12:36:16 CEST----
% \vspace{24pt}
% $$
% \begin{diagram}[w=6ex,h=6ex,tight]
%   I&\topglob{g}\lift{0}{\cel{T}}\dropglob{f}&J\\
%   \dTo<{I}&\lift{-12}{\comm}&\dTo>{J}\\
%   I&\dropglob{f}&J
% \end{diagram}
% \qquad = \qquad
% \begin{diagram}[w=6ex,h=6ex,tight]
%   I&\topglob{g}&J\\
%   \dTo<{I}&\lift{12}{\comm}&\dTo>{J}\\
%   I&\topglob{g}\lift{0}{\cel{T}}\dropglob{f}&J
% \end{diagram}
% \vspace{24pt}
% $$
% %-------------------------------------end of output------checksum:0x2AE3FC99----
with the promised equation \eqref{TXP}.
% $$
% \begin{diagram}[w=6ex,h=6ex,tight]
%   IX&\rTo^{\lambda}&X\\
%   \leftglob{fX}\lift{0}{\cel{T}X}\rightglob{gX}&\shuft{20}{\cel{V}\loft}&\rightglob{X}\\
%   JX&\rTo_{\ell}&X
% \end{diagram}
% \qquad\qquad = \qquad\qquad
% \begin{diagram}[w=6ex,h=6ex,tight]
%   IX&\rTo^{\lambda}&X\\
%   \leftglob{fX}&\shuft{-20}{\cel{U}\loft}&\leftglob{X}\lift{0}{\comm}\rightglob{X}\\
%   JX&\rTo_{\ell}&X
% \end{diagram}
% $$
Since $\cel T$ is cancellable, we can cancel it away to finish.

(iii) implies (i): the arguments in (ii)\impl(iii) can be
reverted: start with \eqref{TXP}, tensor with $\cel T$ on the left, and apply
\eqref{P} to arrive at the axiom for being a $2$-morphism of pseudo-idempotents.
Since both $\cel{U}\loft_X$ and $\cel{V}\loft_X$ are invertible, so is $\cel T X$.
Now take $X$ to be a right cancellable object, and cancel it away to conclude
that already $\cel T$ is invertible.

  Finally, the equivalence 
  (i)\impl(ii')\impl(iii')\impl(i) is completely analogous.
\end{dem}

\begin{cor}\label{unique-2-cell}
  Given two parallel morphisms of units, there is a unique unit $2$-morphism
  between them.
\end{cor}

\begin{dem}
  The  $2$-cell is determined by the previous lemma.
%   Reversing the arguments
%   of the proof shows that that this $2$-cell is  a semi-monoid transformation.
\end{dem}

Next we aim at proving that there is a unit morphism between any two
units.  The strategy is to use the left and right constraints to
produce a unit morphism
$$
I \rTo IJ \rTo J.
$$
As a first step towards this goal we have:
\begin{lemma}\label{IJunit}
  Let $I$ and $J$ be units, and pick a left constraint $\lambda$ for
  $I$ and a right constraint $r$ for $J$.  Put 
  $$
  \gamma \df r_I \lambda_J : IJIJ \to IJ
  $$
Then $(IJ, \gamma)$ is a unit.
\end{lemma}

\begin{dem}
  Since $I$ and $J$ are cancellable, clearly $IJ$ is cancellable too.
  Since $\lambda_J$ and $r_I$ are equi-arrows, $\gamma$ is too.
\end{dem}

\begin{lemma}\label{lambdasemimonoidmap}
  There is a $2$-cell
  \begin{diagram}[w=6ex,h=4.5ex,tight]
  IJIJ & \rTo^{\lambda_J \lambda_J}  & JJ  \\
  \dTo<\gamma  &  \cel{Z}  & \dTo>\beta  \\
  IJ  & \rTo_{\lambda_J}  & J .
  \end{diagram}
Hence $(\lambda_J,\cel{Z})$ is a unit map.  (And there is another 
$2$-cell making $r_I$ a unit map.)
\end{lemma}

\begin{dem}
  The $2$-cell $\cel{Z}$ is defined like this:
  \begin{diagram}[w=48pt,h=36pt,tight,hug]
    IJIJ & &  \\
    \dTo<{IJ\lambda_J}  & \diagonaltopglob{\lambda_J \lambda_J} &  \\
    IJJ & \thintopglob{\lambda_J J}\lift{-2}{\cel{K}^\lambda}\thindropglob{\lambda_{JJ}}
  & JJ  \\
    \thinleftglob{r_I J}\lift{0}{\cel R\inv}\thinrightglob{I\beta}  & 
    \lambda_\beta  & \dTo>\beta  \\
  IJ  & \rTo_{\lambda_J}  & J
  \end{diagram}
  where the $2$-cell $\cel K^\lambda$ is constructed in Lemma~\ref{lambdaXY}.
%   \begin{diagram}[w=6ex,h=4.5ex,tight,hug]
%     IJIJ & &  \\
%     \dTo<{IJ\lambda_J}  &  \rdTo^{\lambda_J \lambda_J} &  \\
%     IJJ & \rTo^{\lambda_J J}  & JJ  \\
%     \dTo<{r_I J}  &  \cel{M}  & \dTo>\beta  \\
%   IJ  & \rTo_{\lambda_J}  & J
%   \end{diagram}
% where $\cel{M}$ is defined from the
% naturality square
% \begin{diagram}[w=6ex,h=4.5ex,tight]
%   IJJ & \rTo^{\lambda_{JJ}}  & JJ  \\
%   \dTo<{I \beta}  &  \cel{N}  & \dTo>\beta  \\
% IJ  & \rTo_{\lambda_J}  & J
% \end{diagram}
% by deforming the left-hand edge using the right modulator for $r$,
% and deforming the top edge using Lemma~\ref{lambdaXY}.
\end{dem}

\begin{cor}\label{exist}
  Given two units, there exists a unit morphism between them.
\end{cor}

\begin{dem}
  Continuing the notation from above,
  by Lemma~\ref{IJunit}, $(IJ,\gamma)$ is a unit, and by 
  Lemma~\ref{lambdasemimonoidmap}, $\lambda: IJ \to J$ is a morphism of
  units.
  Similarly, $r: IJ \to I$ is a unit morphism, and by Lemma~\ref{inv} any
  chosen pseudo-inverse $r\inv : I \to IJ$ is again a unit morphism.
  Finally we take
  $$
  I \stackrel{r\inv}\rTo IJ \stackrel{\lambda}\rTo J.
  $$
\end{dem}

% 
% \begin{BM}
%   Certainly the proof needs to use the left and right actions, since 
%   these are the only means of providing an arrow between different 
%   objects.  It is conceivable that the proof could avoid relying on 
%   the unit structure on the intermediate object $IJ$, but it is 
%   unlikely that it would be any shorter.  In any case, the unit structure
%   on $IJ$ is on independent interest.
% \end{BM}

\begin{satz}[Contractibility]\label{thmC}
  The $2$-category of units in $\CC$ is contractible, if non-empty.
In other words, between any two units there exists a unit morphism, and
between any two parallel unit morphisms there is a unique unit
$2$-morphism.
\end{satz}

\begin{dem}
  By Lemma~\ref{exist} there is a unit morphism between any two units
  (an equi-arrow by definition), and by Lemma~\ref{unique-2-cell}
  there is a unique unit $2$-morphism between any two parallel unit morphisms.
\end{dem}

%%%%%%%%%%%%%%%%%%%%%%%%%%%%%%%%%%%%%%%%%%%%%%%%%%
\section{Classical units}
%%%%%%%%%%%%%%%%%%%%%%%%%%%%%%%%%%%%%%%%%%%%%%%%%%

\label{sec:classical}

In this section we review the classical theory of units
in a monoidal $2$-category, as extracted from the definition
of tricategory of Gordon, Power, and Street~\cite{Gordon-Power-Street}.
In the next section we compare this notion with the cancellable-idempotent
approach of this work.  The equivalence is stated explicitly in 
Theorem~\ref{thmE}.

\begin{blanko}{Tricategories.}
  The notion of tricategory introduced by Gordon, Power, and
  Street~\cite{Gordon-Power-Street} is is roughly a weak category structure
  enriched over bicategories: this means that the structure maps (composition
  and unit) are weak $2$-functors satisfying weak versions of associativity and
  unit constraints.  For the associativity, the pentagon equation is replaced by
  a specified pentagon $3$-cell (TD7), required to satisfy an equation
  corresponding to the $3$-dimensional associahedron.  This equation (TA1) is
  called the nonabelian $4$-cocycle condition.  For the unit structure, three
  families of $3$-cells are specified (TD8): one corresponding to the Kelly
  axiom, one left variant, and one right variant (those two being the
  higher-dimensional analogues of Axioms (2) and (3) of monoidal category).  Two
  axioms are imposed on these three families of $3$-cells: one (TA2) relating
  the left family with the middle family, and one (TA3) relating the right
  family with the middle family.  These are called left and right normalisation.
  (These two axioms are the higher-dimensional analogues of the first argument
  in Kelly's lemma~\ref{Kelly-lemma}.)  It is pointed out in
  \cite{Gordon-Power-Street} that the middle family together with the axioms
  (TA2) and (TA3) completely determine the left and right families if they
  exist.
%   Hence the higher dimensional analogue of
%   Kelly's lemma is not that these cells can be constructed from the middle cell,
%   but rather that if they exist they are unique.
\end{blanko}
  
\begin{blanko}{Monoidal $2$-categories.}
  By specialising the definition of tricategory to the one-object case, and
  requiring everything strict except the units, we arrive at the following
  notion of monoidal $2$-category: 
  a {\em monoidal $2$-category} is a semi-monoidal $2$-category 
  (cf.~\ref{semimonoidal}) equipped with an object $I$,
two natural transformations $\lambda$ and $\rho$
with equi-arrow components
$$
\lambda_X: IX \to X 
$$
$$
\rho_X : XI \to X
$$
and (invertible) $2$-cell data
$$
    \begin{diagram}
    I  X  & \rTo^{\lambda_X}  & X  \\
    \dTo<{I f}  &  \lambda_f  & \dTo>f  \\
    I  Y   & \rTo_{\lambda_Y}  & Y
    \end{diagram}
\qquad\qquad
    \begin{diagram}
    XI  & \rTo^{\rho_X}  & X  \\
    \dTo<{f I}  &  \rho_f  & \dTo>f  \\
    YI   & \rTo_{\rho_Y}  & Y ,
    \end{diagram}
$$
together with three natural modifications $\cel K$, $\cel K^\lambda$,
and $\cel K^\rho$, with invertible components
\begin{align*}
\cel K \; : X \lambda_Y &\Rightarrow \rho_X Y \\
\cel K^\lambda : \lambda_{XY} &\Rightarrow \lambda_X Y \\
\cel K^\rho : X \rho_Y &\Rightarrow \rho_{XY}  .
\end{align*}
We call $\cel K$ the {\em Kelly cell}.
% The three families are denoted (TD8$\mu$), (TD8$\lambda$), and (TD8$\rho$)
% in \cite{Gordon-Power-Street}.

These three families are
subject to the following two equations:
%%  generate triangle \
%%     objects {X \lambda_{YZ}} {\rho_X YZ}  {X\lambda_Y Z} \
%%     arrows {X\cel K^\lambda_{Y,Z}} {\cel K_{X,YZ}} {\cel K_{X,Y} Z}  \
%%     label {TA2}
%-------------------------------------output written 2009/07/12 23:06:10 CEST----
\begin{equation}\label{TA2}
  \xymatrixrowsep{50pt}
  \xymatrixcolsep{42pt}
  \xymatrix @!=0pt {
  X \lambda_{YZ} \ar@{=>}[rr]^{X\cel K^\lambda_{Y,Z}} \ar@{=>}[dr]_{\cel K_{X,YZ}} && X\lambda_Y Z \ar@{=>}[dl]^{\cel K_{X,Y} Z} \\
  & \rho_X YZ &
  }
\end{equation}
%-------------------------------------end of output------checksum:0x81AB62E4----
%%  generate triangle \
%%     objects {X\rho_Y Z} {XY\lambda_Z} {\rho_{XY} Z}  \
%%     arrows {\cel K^\rho_{X,Y} Z} {X\cel K_{Y,Z}} {\cel K_{XY,Z}} \
%%     label {TA3}
%-------------------------------------output written 2009/07/12 23:06:20 CEST----
\begin{equation}\label{TA3}
  \xymatrixrowsep{50pt}
  \xymatrixcolsep{42pt}
  \xymatrix @!=0pt {
  X\rho_Y Z \ar@{=>}[rr]^{\cel K^\rho_{X,Y} Z} \ar@{<=}[dr]_{X\cel K_{Y,Z}} &&
  \rho_{XY} Z \ar@{<=}[dl]^{\cel K_{XY,Z}} \\
  & XY\lambda_Z &
  }
\end{equation}
%-------------------------------------end of output------checksum:0x3F0410A2----
% In other words:
% \vspace{30pt}
% \begin{equation}
%   \label{XIYZ}
% \begin{diagram}[w=8ex,h=6ex,tight]
%   XIYZ&\topglob{X\lambda_Y Z }\lift{-4}{X\cel 
%   K^\lambda_{Y,Z}}\dropglob{X\lambda_{YZ}}&XYZ
% \end{diagram}
% \qquad = \qquad
% \begin{diagram}[w=10ex,h=6ex,tight]
%   XIYZ&\bigtopglob{X\lambda_Y Z}\lift{15}{\cel 
%   K\inv_{X,Y}Z}\lift{-15}{\cel K_{X,YZ}}\straighthor{\rho_X YZ}\bigdropglob{X\lambda_{YZ}}&XYZ
% \end{diagram}
% \vspace{30pt}
% \end{equation}
%   and there is a similar equation involving $\cel K^\rho$.
% (These equations are denoted (TA2) and (TA3) in \cite{Gordon-Power-Street},
% and are called normalised $4$-cocycle conditions.)
\end{blanko}

\begin{BM}
  We have made one change compared to \cite{Gordon-Power-Street}, namely the
  direction of the arrow $\rho_X$: from the viewpoint of $\alpha$ it seems more
  practical to work with $\rho_X:XI\to X$ rather than with the convention of
  $\rho_X : X \to XI$ chosen in \cite{Gordon-Power-Street}.  Since in any case
  it is an equi-arrow, the difference is not essential.  (Gurski in his
  thesis~\cite{Gurski:PhD} has studied a version of tricategory where all the
  equi-arrows in the definition are equipped with specified pseudo-inverses.
  This has the advantage that the definition becomes completely algebraic, in a
  technical sense.)
\end{BM}

\begin{lemma}\label{Icanc}
  The object $I$ is cancellable 
  (independently of the existence of $\cel K$, $\cel K^\lambda$, and $\cel 
  K^\rho$.)
\end{lemma}

\begin{dem}
  We need to establish that `tensoring with $I$ on the left',
  $$
  \mathbb{L} : \Hom(X,Y) \to \Hom(IX,IY) ,
  $$
  is an equivalence of categories.  But this follows since the diagram
  \begin{diagram}[w=9ex,h=6ex,tight]
  \Hom(X,Y) & \rTo^{\mathbb{L}}  & \Hom(IX,IY)  \\
  \dTo<\Id  &    & \dTo>{ \_ \;\#\; \lambda_Y}  \\
  \Hom(X,Y)  & \rTo_{\lambda_X \;\#\; \_}  & \Hom(IX,Y)
  \end{diagram}
  is commutative up to isomorphism: the component at $f:X\to Y$ of this
  isomorphism is just the naturality square $\lambda_f$.  Since the functors
  $\lambda_X \,\#\, \_$ and $ \_ \,\#\, \lambda_Y$ are equivalences, it 
  follows from this isomorphism that $\mathbb{L}$ is too.
\end{dem}

\begin{blanko}{Coherence of the Kelly cell.}\label{Klambda-from-K}
  As remarked in~\cite{Gordon-Power-Street}, if the $\cel K^\lambda$ and $\cel
  K^\rho$ exist, they are determined uniquely from $\cel K$ and the two axioms.
  Indeed, the two equations
%%  generate twotriangles \
%%     objects {I \lambda_{YZ}} {\rho_I YZ}  {I\lambda_Y Z} \
%%     arrows {I\cel K^\lambda_{Y,Z}} {\cel K_{I,YZ}} {\cel K_{I,Y} Z}  \
%%     objects {X\rho_Y I} {XY\lambda_I} {\rho_{XY} I}  \
%%     arrows {\cel K^\rho_{X,Y} I} {X\cel K_{Y,I}} {\cel K_{XY,I}} \
%%     label {KKK}
%-------------------------------------output written 2009/07/15 14:25:19 CEST----
\begin{equation}\label{KKK}
  \xymatrixrowsep{50pt}
  \xymatrixcolsep{42pt}
  \xymatrix @!=0pt {
  I \lambda_{YZ} \ar@{=>}[rr]^{I\cel K^\lambda_{Y,Z}} \ar@{=>}[dr]_{\cel K_{I,YZ}} && I\lambda_Y Z \ar@{=>}[dl]^{\cel K_{I,Y} Z} \\
  & \rho_I YZ &
  }
  \qquad\qquad
   \xymatrixrowsep{50pt}
  \xymatrixcolsep{42pt}
  \xymatrix @!=0pt {
  X\rho_Y I \ar@{=>}[rr]^{\cel K^\rho_{X,Y} I} \ar@{=>}[dr]_{X\cel K_{Y,I}} && \rho_{XY} I \ar@{=>}[dl]^{\cel K_{XY,I}} \\
  & XY\lambda_I &
  }
\end{equation}
%-------------------------------------end of output------checksum:0x3A21A935----
which are just special cases of \eqref{TA2} and \eqref{TA3} uniquely determine 
$\cel K^\lambda$ and $\cel K^\rho$, by cancellability of $I$.
But these two special cases of the axioms do not imply the general case.

We shall take the Kelly cell $\cel K$ as the main structure, and say that $\cel
K$ is {\em coherent on the left} (resp.~{\em on the right}) if Axiom~\eqref{TA2}
(resp.~\eqref{TA3}) holds for the induced cell $\cel K^\lambda$ (resp.~$\cel 
K^{\rho}$).  We just say {\em coherent} if both hold.  We shall
see (\ref{KLR}) that in fact coherence on the left implies coherence on the
right and vice versa.
\end{blanko}

\begin{blanko}{Naturality.}\label{natK}
  The Kelly cell is a modification.  For future reference we spell out the
  naturality condition satisfied: given arrows 
  $f: X \to X'$ and $g: Y \to Y'$, we have
%%  generate verticalCylinder \
%%     startobjects XIY XY \
%%     endobjects X'IY' X'Y' \
%%     straightarrows {fIg} {fg} \
%%     curvedarrows {\rho_X Y} {X\lambda_Y} {\rho_{X'}Y'} {X'\lambda_{Y'}} \
%%     twocells {\cel K_{X,Y}} {\cel K_{X',Y'}} {f \lambda_g} {\rho_f g} \
%%     label {}
%-------------------------------------output written 2009/07/15 14:31:32 CEST----
\vspace{24pt}
$$
\begin{diagram}[w=6ex,h=6ex,tight]
  XIY&\topglob{\rho_X Y}\lift{0}{\cel K_{X,Y}}\dropglob{X\lambda_Y}&XY\\
  \dTo<{fIg}&\lift{-12}{f \lambda_g}&\dTo>{fg}\\
  X'IY'&\dropglob{X'\lambda_{Y'}}&X'Y'
\end{diagram}
\qquad = \qquad
\begin{diagram}[w=6ex,h=6ex,tight]
  XIY&\topglob{\rho_X Y}&XY\\
  \dTo<{fIg}&\lift{12}{\rho_f g}&\dTo>{fg}\\
  X'IY'&\topglob{\rho_{X'}Y'}\lift{0}{\cel K_{X',Y'}}\dropglob{X'\lambda_{Y'}}&X'Y'
\end{diagram}
\vspace{24pt}
$$
\end{blanko}

\begin{BM}
  Particularly useful is naturality of $\lambda$ with respect to $\lambda_X$
  and naturality of $\rho$ with respect to $\rho_X$.
  In these cases, since $\lambda_X$ and $\rho_X$ are equi-arrows, we can
  cancel them and find the following invertible $2$-cells:
  \begin{align*}
    \cel N^\lambda : I \lambda_X &\Rightarrow \lambda_{IX} \\
    \cel N^\rho\,: \,\rho_{XI} &\Rightarrow X\rho_I ,
  \end{align*}
  in analogy with Observation~(5) of monoidal categories.
\end{BM}

The following lemma holds for $\cel K$ independently of Axioms~\eqref{TA2} and 
\eqref{TA3}:
\begin{lemma}\label{NK=KN}
  The Kelly cell $\cel K$ satisfies the equation
%   We have the following square of invertible $2$-cells:
%   \[
%     \xymatrixrowsep{50pt}
%   \xymatrixcolsep{64pt}
%   \xymatrix @!=0pt {
%   XI\lambda_Y \ar@{=>}[r]^{\cel K_{XI,Y}} \ar@{=>}[d]_{X\cel N^\lambda} & 
%   \rho_{XI} Y 
%   \ar@{=>}[d]^{\cel N^\rho Y}\\
%   X\lambda_{IY} \ar@{=>}[r]_{\cel K_{X,IY}} & \rho_X IY .
%   }\]
% \end{lemma}
% 
% \begin{dem}
%   The equation is
\vspace{30pt}
$$
\begin{diagram}[w=9ex,h=7ex,tight]
  XIIY&\bigtopglob{\rho_X IY}\lift{15}{\cel{K}_{X,IY}}\lift{-12}{X\cel 
  N^\lambda}\straighthor{X\lambda_{IY}}\bigdropglob{XI\lambda_Y}&XIY
\end{diagram}
  \qquad = \qquad
\begin{diagram}[w=9ex,h=7ex,tight]
  XIIY&\bigtopglob{\rho_X IY}\lift{15}{\cel N^\rho Y}\lift{-15}{\cel K_{XI,Y}}
  \straighthor{\rho_{XI}Y}\bigdropglob{XI\lambda_Y}&XIY
\end{diagram}
\vspace{30pt}
$$
\end{lemma}

\begin{dem}
  It is enough to establish this equation after post-whiskering with
  $X\lambda_Y$.  The rest is a routine calculation, using on one side the
  definition of the cell $\cel N^\lambda$, then naturality of $\cel K$ with
  respect to $f=X$ and $g=\lambda_Y$.  On the other side, use the definition of
  $\cel N^\rho$ and then naturality of $\cel K$ with respect to $f=\rho_X$ and
  $g=Y$.  In the end, two $\cel K$-cells cancel.
\end{dem}

Combining the $2$-cells described so far we get
$$
\rho_I I \stackrel{\cel K\inv}{\Rightarrow} I\lambda_I
\stackrel{\cel N^\lambda}{\Rightarrow} \lambda_{II} \stackrel{\cel
K^\lambda}{\Rightarrow} \lambda_I I
$$
and hence, by cancelling $I$ on the right, an invertible $2$-cell
$$
  \cel P : \rho_I \Rightarrow \lambda_I .
$$
Now we could also define $\cel Q : \rho_I \Rightarrow \lambda_I$
in terms of
$$
I\rho_I \stackrel{\cel K^\rho} \Rightarrow
\rho_{II} \stackrel{\cel N^\rho} \Rightarrow
\rho_I I \stackrel{\cel K\inv}\Rightarrow
 I\lambda_I .
$$
Finally,  in analogy with Axiom~(1) for monoidal categories:
\begin{lemma}\label{U=V}
  We have $\cel P = \cel Q$.  (This is true independently of Axioms~\eqref{TA2} 
  and \eqref{TA3}.)
\end{lemma}

\begin{dem}
  Since $I$ is cancellable, it is enough to show $I\cel P I = I\cel Q I$.
  To establish this equation, use the constructions of $\cel P$ and $\cel Q$,
  then substitute the characterising Equations~\eqref{KKK} for the auxiliary cells
  $\cel K^\lambda$ and $\cel K^\rho$, and finally use Lemma~\ref{NK=KN}.
\end{dem}

\begin{blanko}{The $2$-category of GPS units.}\label{GPS}
  For short we shall say {\em GPS unit} for the notion of unit just introduced.
  In summary, a GPS unit is a quadruple $(I, \lambda, \rho, \cel K)$ where $I$
  is an object, $\lambda_X$ and $\rho_X$ are natural transformations
  with equi-arrow components, and $\cel K : X\lambda_Y
  \Rightarrow \rho_X Y$ is a good Kelly cell (natural in $X$ and $Y$, of course).

  A {\em morphism of GPS units} from $(I,\lambda,\rho,\cel K)$ to $(J,\ell, r,
  \cel H)$ is an arrow $u: I \to J$ equipped with natural families of invertible
  $2$-cells
$$
\begin{diagram}
  IX & \rTo^{uX}  & JX  \\
  \dTo<{\lambda_X}  & \cel{U}\loft_X   & \dTo>{\ell_X}  \\
  X  & \rTo_X  & X
  \end{diagram}
  \qquad\qquad
  \begin{diagram}
  XI & \rTo^{Xu}  & XJ  \\
  \dTo<{\rho_X}  & \cel{U}\rught_X   & \dTo>{r_X}  \\
  X  & \rTo_X  & X
  \end{diagram}
$$
satisfying the equation
%%  generate horizontalCylinder \
%%     startobjects XIY XY \
%%     endobjects XJY XY \
%%     straightarrows {XuY} {XY} \
%%     curvedarrows {X\lambda_Y} {\rho_X Y} {X\ell_Y} {r_X Y} \
%%     twocells {\cel K} {\cel H} {\cel{U}\rught_X Y} {X\cel{U}\loft_Y} \
%%     label {PK}
%-------------------------------------output written 2009/07/13 15:05:51 CEST----
\begin{equation}\label{PK}
\begin{diagram}[w=6ex,h=6ex,tight]
  XIY&\rTo^{XuY}&XJY\\
  \leftglob{X\lambda_Y}\lift{0}{\cel K}\rightglob{\rho_X Y}&\shuft{20}{\cel{U}\rught_X Y}&\rightglob{r_X Y}\\
  XY&\rTo_{XY}&XY
\end{diagram}
\qquad\qquad = \qquad\qquad
\begin{diagram}[w=6ex,h=6ex,tight]
  XIY&\rTo^{XuY}&XJY\\
  \leftglob{X\lambda_Y}&\shuft{-20}{X\cel{U}\loft_Y}&\leftglob{X\ell_Y}\lift{0}{\cel H}\rightglob{r_X Y}\\
  XY&\rTo_{XY}&XY
\end{diagram}
\end{equation}
%-------------------------------------end of output------checksum:0x2277B34B----

  Finally, a {\em $2$-morphism of GPS unit maps} is a $2$-cell $\cel T: u
  \Rightarrow v$ satisfying the compatibility conditions \eqref{TXP} and 
  \eqref{TXQ} of Lemma~\ref{unit2map}.
\end{blanko}

\begin{blanko}{Remarks.}
  Note first that $u$ is automatically an equi-arrow.  Observe also that
  $\cel{U}\loft$ and $\cel{U}\rught$ completely determine each other by
  Equation~\eqref{PK}, as is easily seen by taking respectively $X$ to be a left
  cancellable object and $Y$ to be a right cancellable object.  Finally note
  that 
  there are two further equations, expressing compatibility with $\cel K^\lambda$
  and $\cel K^\rho$, but they can be deduced from Equation~\eqref{PK},
  independently of the coherence Axioms~\eqref{TA2} and \eqref{TA3}.
  Here is the one for $\cel K^\lambda$
  for future reference:
%%  generate horizontalCylinder \
%%     startobjects IXY XY \
%%     endobjects JXY XY \
%%     straightarrows {uXY} {XY} \
%%     curvedarrows {\lambda_{XY}} {\lambda_X Y} {\ell_{XY}} {\ell_X Y} \
%%     twocells {\cel K^\lambda} {\cel H^\ell} {\cel{U}\loft_X Y} {\cel{U}\loft_{XY}} \
%%     label {KP=PH}
%-------------------------------------output written 2009/07/13 15:06:17 CEST----
\begin{equation}\label{KP=PH}
\begin{diagram}[w=6ex,h=6ex,tight]
  IXY&\rTo^{uXY}&JXY\\
  \leftglob{\lambda_{XY}}\lift{0}{\cel K^\lambda}\rightglob{\lambda_X Y}&\shuft{20}{\cel{U}\loft_X Y}&\rightglob{\ell_X Y}\\
  XY&\rTo_{XY}&XY
\end{diagram}
\qquad\qquad = \qquad\qquad
\begin{diagram}[w=6ex,h=6ex,tight]
  IXY&\rTo^{uXY}&JXY\\
  \leftglob{\lambda_{XY}}&\shuft{-20}{\cel{U}\loft_{XY}}&\leftglob{\ell_{XY}}\lift{0}{\cel H^\ell}\rightglob{\ell_X Y}\\
  XY&\rTo_{XY}&XY
\end{diagram}
\end{equation}
%-------------------------------------end of output------checksum:0x43BA2E78----
\end{blanko}

%%%%%%%%%%%%%%%%%%%%%%%%%%%%%%%%%%%%%%%%%%%%%%%%%%
\section{Comparison with classical theory (Theorem~\ref{thmE})}
%%%%%%%%%%%%%%%%%%%%%%%%%%%%%%%%%%%%%%%%%%%%%%%%%%
\label{sec:comparison}

In this section we prove the equivalence between the two notions of unit.

\begin{blanko}{From cancellable-idempotent units to GPS units.}
  We fix a unit object $(I,\alpha)$.  We also
  assume chosen a left constraint $\lambda_X : IX \to X$ with $\cel L_X:
  I\lambda_X \Rightarrow \alpha X$, and a right constraint $\rho_X : XI \to X$
  with $\cel R_X: X \alpha \Rightarrow \rho_X I$.  First of all, in analogy with
   Axioms (2) and (3) of monoidal categories we have:
\end{blanko}

\begin{lemma}\label{lambdaXY}
  There are unique natural invertible $2$-cells
  \begin{align*}
    \cel{K}^\lambda : \lambda_{XY} &\Rightarrow \lambda_X Y \\
  \cel K^\rho:  X \rho_Y &\Rightarrow \rho_{XY}
  \end{align*}
  satisfying
\vspace{30pt}
\begin{equation}
  \label{IIAB}
\begin{diagram}[w=6ex,h=6ex,tight]
  IIXY&\topglob{I\lambda_X Y}\lift{0}{I\cel{K}^\lambda}\dropglob{I\lambda_{XY}}&IXY
\end{diagram}
\qquad = \qquad
\begin{diagram}[w=8ex,h=6ex,tight]
  IIXY&\bigtopglob{I\lambda_X Y}\lift{12}{\cel{L}\inv Y}\lift{-12}{\cel{L}}\straighthor{\alpha XY}\bigdropglob{I\lambda_{XY}}&IXY
\end{diagram}
\vspace{30pt}
\end{equation}

\bigskip
\bigskip

\begin{equation}
  \label{XYII}
\begin{diagram}[w=6ex,h=6ex,tight]
  XYII&\topglob{\rho_{XY}I}\lift{0}{\cel K^\rho I}\dropglob{X\rho_Y I}&XYI
\end{diagram}
\qquad = \qquad
\begin{diagram}[w=8ex,h=6ex,tight]
  XYII&\bigtopglob{\rho_{XY}I}\lift{12}{\cel R}\lift{-12}{X\cel R\inv}\straighthor{XY\alpha}\bigdropglob{X\rho_Y I}&XYI
\end{diagram}
\vspace{30pt}
\end{equation}
\end{lemma}

\begin{dem}
  The condition precisely defines the $2$-cell, since $I$ is
  cancellable.  
\end{dem}

% \footnotesize
% Remark: the operation of cancelling away an $I$-factor for the purpose
% of defining a $2$-cell depends crucially on naturality, and in fact
% the
% constructed $2$-cell involves $\Lambda$ in ites definition although it
% is not apparent in the argument above.  To spell it out goes like
% this: consider the two naturality squares:
% $$
% \begin{diagram}[w=6ex,h=4.5ex,tight]
% IIXY & \rTo^{\lambda_{IXY}}  & IXY  \\
% \dTo<{I \lambda_{XY}}  &  \Lambda_{XY}  & \dTo>{\lambda_{XY}}  \\
% IXY  & \rTo_{\lambda_{XY}}  & XY
% \end{diagram}
% \qquad
% \text{ and }
% \qquad
% \begin{diagram}[w=6ex,h=4.5ex,tight]
% IIXY & \rTo^{\lambda_{IXY}}  & IXY  \\
% \dTo<{I \lambda_{X}Y}  &  \Lambda_X Y  & \dTo>{\lambda_{X}Y}  \\
% IXY  & \rTo_{\lambda_{XY}}  & XY
% \end{diagram}
% $$
% These can be glued together along top and bottom, and at the left-hand
% end of the cylinder we can put the $2$-cell of (\ref{IIAB}).  This
% defines the right-hand $2$-cell.
% \normalsize

% \begin{lemma}
%   There are natural invertible $2$-cells
%   $$
%   \lambda_{I X} \Leftrightarrow  I
%    \lambda_X
%   $$
%   $$
%   \rho_{X  I} \Leftrightarrow \rho_X
%    I
%   $$
% \end{lemma}
% \begin{dem}
%   This follows from  the naturality of the left action $\lambda$
%   corresponding to the arrow $\lambda_X$, and from the naturality of 
%   the right action $\rho$ with respect to $\rho_X$.
% %   \begin{diagram}[w=6ex,h=4.5ex,tight]
% %   IIX & \rTo^{\lambda_{I X}}  & IX  \\
% %   \dTo<{I \lambda_X}  &    & \dTo>{\lambda_X}  \\
% %   IX  & \rTo_{\lambda_X}  & X
% %   \end{diagram}
% %   by cancelling away $\lambda_X$.
% \end{dem}

\begin{lemma}
  For fixed left constraint $(\lambda,\cel L)$ and fixed right
  constraint $(\rho,\cel R)$, 
  there is a canonical family of invertible $2$-cells  (the Kelly cell)
  $$
  \cel K : X  \lambda_Y \Rightarrow \rho_X  Y ,
  $$
  natural in $X$ and $Y$.
\end{lemma}

\begin{dem}
  This is analogous to the construction of the associator:
  $\cel K$ is defined as the unique $2$-cell $\cel K: X\lambda_Y
  \Rightarrow \rho_X Y$ satisfying the equation
%   
%   start with the commutative square
%   \begin{equation}\label{KellyC} 
%   \begin{diagram}[w=6ex,h=4.5ex,tight]
%   XIIY & \rTo^{XI\lambda_Y}  & XIY  \\
%   \dTo<{\rho_X IY}  &  \comm  & \dTo>{\rho_X Y}  \\
%   XIY  & \rTo_{X\lambda_Y}  & XY
%   \end{diagram}
% \end{equation}
%   and attach cells as follows:
  \vspace{24pt}
  \begin{equation}\label{RYXL}
  \begin{diagram}[w=52pt,h=39pt,tight]
  XIIY & \thintopglob{X\alpha
  Y}\lift{0}{X\cel{L}}\thindropglob{XI\lambda_Y}  & 
  XIY  \\
  \thinleftglob{X\alpha Y}\lift{0}{\cel{R}Y}\thinrightglob{\rho_X IY}  &  
  \comm
  & \dTo>{\rho_X Y}  \\
  XIY  & \rTo_{X\lambda_Y}  & XY
  \end{diagram}
  \qquad = \qquad
    \begin{diagram}[w=30pt,h=22.5pt,tight,hug]
  XIIY &&&&  \\
  &\rdTo^{X\alpha Y} &&&& \\
  && XIY && \\
  &&&\diagonaltopglob{\rho_X Y}
  \diagonaldropglob{X\lambda_Y}\NWSElabel{\cel K}   & \\
  &&&& XY
  \end{diagram}
  \vspace{24pt}
\end{equation}
   This makes sense    since $X\alpha Y$ is an equi-arrow, so we can cancel it 
   away.  Clearly $\cel K$ is invertible since $\cel L$ and $\cel R$ are.
\end{dem}

We constructed $\cel K^\lambda$ and $\cel K^\rho$ directly from $\cel L$,
and $\cel R$.  Meanwhile we also constructed $\cel K$,
and we know from classical theory (\ref{Klambda-from-K}) that this cell determines the two 
others.  The following proposition shows that all these constructions match up,
and in particular that the constructed Kelly cell is coherent on both sides:
% 
% As a corollary to the proposition, here is a characterisation in terms
% of $\cel K$:
% \begin{cor}\label{WKK}
%   We have commutative triangles of $2$-cells:
% %%  generate twotriangles \
% %%     objects {I \lambda_{XY}} {\rho_I XY} {I\lambda_X Y} \
% %%     arrows  {I \cel K^\lambda} {\cel K} {\cel K Y}  \
% %%     objects {X\rho_Y I} {XY\lambda_I} {\rho_{XY}I} \
% %%     arrows  {\cel K^\rho I} {X\cel K} {\cel K} \
% %%     label {}
% %-------------------------------------output written 2009/07/08 15:14:47 CEST----
% $$
%   \xymatrixrowsep{50pt}
%   \xymatrixcolsep{42pt}
%   \xymatrix @!=0pt {
%   I \lambda_{XY} \ar@{=>}[rr]^{I \cel K^\lambda} \ar@{=>}[dr]_{\cel K} && I\lambda_X Y \ar@{=>}[dl]^{\cel K Y} \\
%   & \rho_I XY &
%   }
%   \qquad\qquad
%    \xymatrixrowsep{50pt}
%   \xymatrixcolsep{42pt}
%   \xymatrix @!=0pt {
%   X\rho_Y I \ar@{=>}[rr]^{\cel K^\rho I} \ar@{<=}[dr]_{X\cel K} && \rho_{XY}I 
%   \ar@{<=}[dl]^{\cel K} \\
%   & XY\lambda_I &
%   }
% $$
% %-------------------------------------end of output------checksum:0x44CA54A7----
% \end{cor}
% \begin{dem}
%   This follows directly by combining \ref{KIW} with (\ref{IIAB}).
% \end{dem}
% 
% 
% 
% The Corollary is just a special case of the next, more difficult, result:

\begin{prop}\label{Kcoh}
  The families of $2$-cells constructed, $\cel K$, $\cel K^\lambda$ and $\cel 
  K^\rho$ satisfy the GPS unit axioms~\eqref{TA2} and \eqref{TA3}:
%%  generate twotriangles \
%%     objects {X \lambda_{YZ}} {\rho_X YZ}  {X\lambda_Y Z} \
%%     arrows {X\cel K^\lambda_{Y,Z}} {\cel K_{X,YZ}} {\cel K_{X,Y} Z}  \
%%     objects {X\rho_Y Z} {XY\lambda_Z} {\rho_{XY} Z}  \
%%     arrows {\cel K^\rho_{X,Y} Z} {X\cel K_{Y,Z}} {\cel K_{XY,Z}} \
%%     label {}
%-------------------------------------output written 2009/07/08 15:17:16 CEST----
$$
  \xymatrixrowsep{50pt}
  \xymatrixcolsep{42pt}
  \xymatrix @!=0pt {
  X \lambda_{YZ} \ar@{=>}[rr]^{X\cel K^\lambda_{Y,Z}} \ar@{=>}[dr]_{\cel K_{X,YZ}} && X\lambda_Y Z \ar@{=>}[dl]^{\cel K_{X,Y} Z} \\
  & \rho_X YZ &
  }
  \qquad\qquad
   \xymatrixrowsep{50pt}
  \xymatrixcolsep{42pt}
  \xymatrix @!=0pt {
  X\rho_Y Z \ar@{=>}[rr]^{\cel K^\rho_{X,Y} Z} \ar@{<=}[dr]_{X\cel K_{Y,Z}} &&
  \rho_{XY} Z \ar@{<=}[dl]^{\cel K_{XY,Z}} \\
  & XY\lambda_Z &
  }
$$
%-------------------------------------end of output------checksum:0xA62727C7----
\end{prop}

\begin{dem}
  We treat the left constraint (the right constraint being completely 
  analogous).
  We need to establish
  \vspace{30pt}
$$
\begin{diagram}[w=10ex,h=6ex,tight]
  XIYZ&\bigtopglob{\rho_X YZ}\lift{15}{\cel K_{X,Y}Z}
  \lift{-17}{X\cel 
  K^\lambda_{Y,Z}}\straighthor{X\lambda_Y Z}\bigdropglob{X\lambda_{YZ}}&XYZ
\end{diagram}
\qquad = \qquad
\begin{diagram}[w=8ex,h=6ex,tight]
  XIYZ&\topglob{\rho_X YZ}\lift{-2}{\cel K_{X,YZ}}\dropglob{X\lambda_{YZ}}&XYZ
\end{diagram}
\vspace{30pt}
$$
and it is enough to establish this equation pre-whiskered with $X\alpha Y Z$.
In the diagram resulting from the left-hand side:
  \vspace{30pt}
$$
\begin{diagram}[w=10ex,h=6ex,tight]
  XIIYZ & \rTo^{X\alpha YZ} &
  XIYZ&\bigtopglob{\rho_X YZ}\lift{15}{\cel K_{X,Y}Z}\lift{-17}{X\cel 
  K^\lambda_{Y,Z}}\straighthor{X\lambda_Y Z}\bigdropglob{X\lambda_{YZ}}&XYZ
\end{diagram}
\vspace{30pt}
$$
we can replace $(X\alpha YZ) \# (\cel K_{X,Y}Z)$ by the expression that
defined  $\cel 
K_{X,Y}Z$ (cf.~\ref{RYXL}), yielding altogether
\begin{diagram}[w=48pt,h=36pt,tight,scriptlabels,hug]
 &   & XIYZ & & \\
  & \NEdiagonaltopglob{X \alpha YZ} 
  \NEdiagonaldropglob{XI\lambda_Y Z}\NESWlabel{X \cel L_Y Z}   &  & \rdTo^{\rho_X YZ}  &  \\
XIIYZ  &    & \comm & & XYZ \\
  & \diagonaltopglob{\rho_X IYZ}
  \diagonaldropglob{X\alpha YZ}\NWSElabel{\cel{R}_X YZ}  & &\NEdiagonaltopglob{X 
  \lambda_Y Z} 
  \NEdiagonaldropglob{X\lambda_{YZ}}\NESWlabel{X \cel K^\lambda_{Y,Z}} & \\
  && XIYZ
\end{diagram}
Here we can move the cell $X\cel K^\lambda_{Y,Z}$ across the square,
where it becomes $XI \cel K^\lambda_{Y,Z}$ and combines with $X\cel L_Y Z$
to give altogether $X \cel L_{YZ}$ (cf.~\eqref{IIAB}).
The resulting diagram
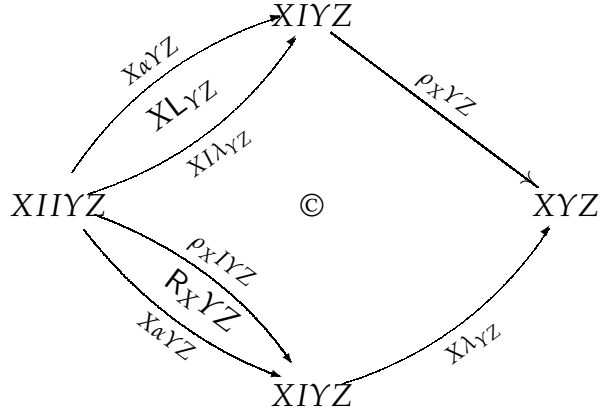
\begin{diagram}[w=48pt,h=36pt,tight,scriptlabels,hug]
 &   & XIYZ & & \\
  & \NEdiagonaltopglob{X \alpha YZ} 
  \NEdiagonaldropglob{XI\lambda_{YZ}}\NESWlabel{X \cel L_{YZ}}   &  & \rdTo^{\rho_X YZ}  &  \\
XIIYZ  &    & \comm & & XYZ \\
  & \diagonaltopglob{\rho_X IYZ}
  \diagonaldropglob{X\alpha YZ}\NWSElabel{\cel{R}_X YZ}  & &
  \NEdiagonaldropglob{X\lambda_{YZ}} & \\
  && XIYZ
\end{diagram}
is nothing but
  
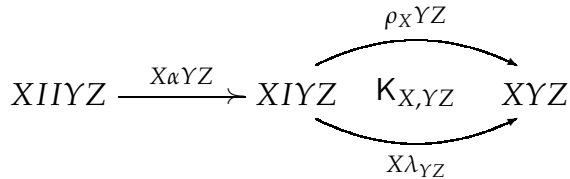
\vspace{30pt}
$$
\begin{diagram}[w=8ex,h=6ex,tight]
  XIIYZ & \rTo^{X\alpha YZ} &
  XIYZ&\topglob{\rho_X YZ}\lift{-2}{\cel 
  K_{X,YZ}}\dropglob{X\lambda_{YZ}}&XYZ
\end{diagram}
\vspace{30pt}
$$
(by Equation~\eqref{RYXL} again)
which is what we wanted to establish.
\end{dem}

Hereby we have concluded the construction of a GPS unit from $(I,\alpha)$.
We will also need a result for morphisms:

\begin{prop}\label{Kelly-compatibility}
  Let $(u,\cel U) :(I,\alpha) \to (J,\beta)$ be a morphism of units in the sense
  of \ref{main-def}, and consider the two canonical $2$-cells $\cel{U}\loft$ and
  $\cel{U}\rught$ constructed in Lemma~\ref{unitmap}.  Then Equation~\eqref{PK}
  holds:
%%  generate horizontalCylinder \
%%     startobjects XIY XY \
%%     endobjects XJY XY \
%%     straightarrows {XuY} {XY} \
%%     curvedarrows {X\lambda_Y} {\rho_X Y} {X\ell_Y} {r_X Y} \
%%     twocells {\cel K} {\cel H} {\cel{U}\rught_X Y} {X\cel{U}\loft_Y} \
%%     label {}
%-------------------------------------output written 2009/07/15 14:38:49 CEST----
$$
\begin{diagram}[w=6ex,h=6ex,tight]
  XIY&\rTo^{XuY}&XJY\\
  \leftglob{X\lambda_Y}\lift{0}{\cel K}\rightglob{\rho_X Y}&\shuft{20}{\cel{U}\rught_X Y}&\rightglob{r_X Y}\\
  XY&\rTo_{XY}&XY
\end{diagram}
\qquad\qquad = \qquad\qquad
\begin{diagram}[w=6ex,h=6ex,tight]
  XIY&\rTo^{XuY}&XJY\\
  \leftglob{X\lambda_Y}&\shuft{-20}{X\cel{U}\loft_Y}&\leftglob{X\ell_Y}\lift{0}{\cel H}\rightglob{r_X Y}\\
  XY&\rTo_{XY}&XY
\end{diagram}
$$
%-------------------------------------end of output------checksum:0xCE7142EF----
  (Hence $(u,\cel{U}\loft,\cel{U}\rught)$ is a morphism of GPS units.)
\end{prop}
\begin{dem}
  It is enough to prove the equation obtained by pasting the $2$-cell $X\cel{U}
  Y$ on top of each side of the equation.  This enables us to use the
  characterising equation for $\cel K$ and $\cel H$.  After this rewriting, we
  are in position to apply Equations~\eqref{P} and \eqref{Q},
  and after cancelling $\cel R$ and $\cel L$ cells, the
  resulting equation amounts to a cube, where it is easy to see that each side
  is just $\cel{U}\rught_X \cel{U}\loft_Y$.
\end{dem}

\begin{blanko}{From GPS units to cancellable-idempotent units.}
  Given a GPS unit $(I, \lambda, \rho, \cel K)$, just put 
  $$
  \alpha \df \lambda_I ,
  $$
  then $(I,\alpha)$ is a unit object in the sense of \ref{main-def}.  Indeed, we
  already observed that $I$ is cancellable (\ref{Icanc}), and from the outset
  $\lambda_I$ is an equi-arrow.  That's all!  To construct it we didn't even
  need the Kelly cell, or any of the auxiliary cells or their axioms.
\end{blanko}

\begin{blanko}{Left and right actions from the Kelly cell.}\label{LRfromK}
  Start with natural left and right constraints $\lambda$ and
  $\rho$ and a Kelly cell $\cel K : X\lambda_Y \Rightarrow \rho_X Y$ (not 
  required to be coherent on either side).  Construct $\cel K^\lambda$ as in 
  \ref{Klambda-from-K}, put $\alpha \df \lambda_I$, and define left and right 
  actions as follows.
We define $\cel L_X$ as
$$
I\lambda_X \stackrel{\cel N^\lambda}{\Rightarrow} \lambda_{IX} 
\stackrel{\cel K^\lambda}{\Rightarrow} \lambda_I X = \alpha X ,
$$
while we define $\cel R_X$ simply as
$$
 X\alpha =  X\lambda_I \stackrel{\cel K_{X,I}}{\Rightarrow}   \rho_X I .
$$
\end{blanko}

\begin{prop}\label{KLR}
  For fixed $(I,\lambda,\rho,\cel K)$, the following are equivalent:
  \begin{punkt-i}
    \item The left and right $2$-cells $\cel L$ and $\cel R$ just constructed in
    \ref{LRfromK} are compatible with the Kelly cell in the sense of
    Equation~\eqref{RYXL}.
  
    \item The Kelly cell $\cel K$ is coherent on the left
    (i.e.~satisfies Axiom~\eqref{TA2}).
    
    \item[(ii')] The Kelly cell $\cel K$ is coherent on the right
    (i.e.~satisfies Axiom~\eqref{TA3}).
  \end{punkt-i}
\end{prop}

\begin{dem}
  Proposition~\ref{Kcoh} already says that (i) implies both (ii) and (ii').
  To prove (ii)\impl(i),
  we start with an auxiliary observation: by massaging the naturality equation
%%  generate verticalCylinder \
%%     startobjects XIIY XIY \
%%     endobjects XIY XY \
%%     straightarrows {XI\lambda_Y} {X\lambda_Y} \
%%     curvedarrows {\rho_X IY} {X\lambda_{IY}} {\rho_X Y} {X\lambda_Y} \
%%     twocells {\cel{K}_{X,IY}} {\cel{K}_{X,Y}} {X\lambda_{\lambda_Y}} {\comm} \
%%     label {}
%-------------------------------------output written 2009/07/13 13:35:44 CEST----
\vspace{24pt}
$$
\begin{diagram}[w=6ex,h=6ex,tight]
  XIIY&\topglob{\rho_X IY}\lift{0}{\cel{K}_{X,IY}}\dropglob{X\lambda_{IY}}&XIY\\
  \dTo<{XI\lambda_Y}&\lift{-12}{X\lambda_{\lambda_Y}}&\dTo>{X\lambda_Y}\\
  XIY&\dropglob{X\lambda_Y}&XY
\end{diagram}
\qquad = \qquad
\begin{diagram}[w=6ex,h=6ex,tight]
  XIIY&\topglob{\rho_X IY}&XIY\\
  \dTo<{XI\lambda_Y}&\lift{12}{\comm}&\dTo>{X\lambda_Y}\\
  XIY&\topglob{\rho_X Y}\lift{0}{\cel{K}_{X,Y}}\dropglob{X\lambda_Y}&XY
\end{diagram}
\vspace{24pt}
$$
%-------------------------------------end of output------checksum:0x8CF2F8AF----
we find the equation
\vspace{24pt}
\begin{equation}\label{tailored}
\begin{diagram}[w=7ex,h=7ex,tight]
  XIIY&\topglob{X \lambda_{IY}}\lift{0}{X\cel{N}^\lambda}\dropglob{XI\lambda_Y}&XIY\\
  \dTo<{\rho_X IY}&\lift{-12}{\comm}&\dTo>{\rho_X Y}\\
  XIY&\rTo_{X\lambda_Y}&XY
\end{diagram}
\qquad = \qquad
\begin{diagram}[w=7ex,h=7ex,tight]
  XIIY&\topglob{X \lambda_{IY}}\lift{0}{\cel{K}\inv_{X,IY}}\dropglob{\rho_X IY}&XIY\\
  \dTo<{\rho_X IY}&\lift{-12}{\comm}&
  \leftglob{X\lambda_Y}\lift{0}{\cel K_{X,Y}}\rightglob{\rho_X Y}\\
  XIY&\rTo_{X\lambda_Y}&XY
\end{diagram}
\vspace{24pt}
\end{equation}
tailor-made to a substitution we shall perform in a moment.  

Now for the main computation, assuming first that $\cel K$ is coherent on the
left, i.e.~that Axiom~\eqref{TA2} holds.  Start with the left-hand side of
Equation~\eqref{RYXL}, and insert the definitions we made for $\cel L$ and $\cel
R$ to arrive at 
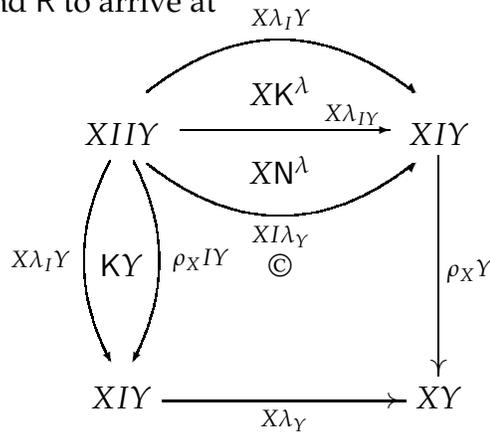
\vspace{32pt}
$$
\begin{diagram}[w=60pt,h=50pt,tight]
  XIIY & \bigtopglob{X\lambda_I Y}
  \lift{15}{X\cel{K}^\lambda} 
  \straighthor{X\lambda_{IY}}
  \lift{-15}{X\cel{N}^\lambda}
  \bigdropglob{XI\lambda_Y}  & 
  XIY  \\
  \thinleftglob{X\lambda_I Y}\lift{0}{\cel{K}Y}\thinrightglob{\rho_X IY}  &  
  \comm
  & \dTo>{\rho_X Y}  \\
  XIY  & \rTo_{X\lambda_Y}  & XY
  \end{diagram}
  \vspace{24pt}
$$
in which we can now substitute \eqref{tailored} to get
  
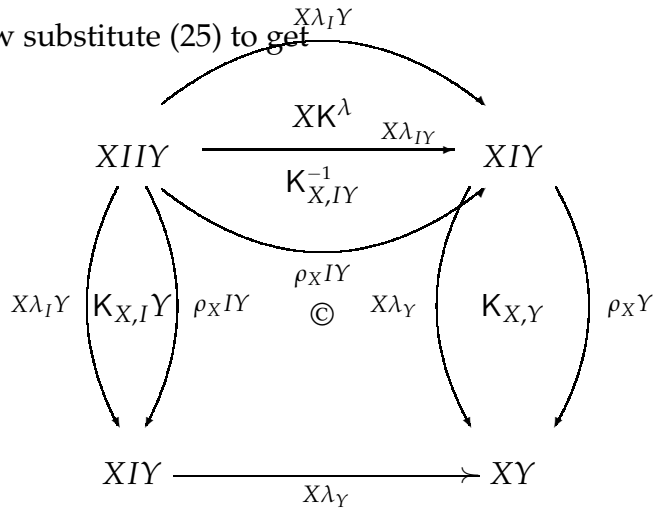
\vspace{24pt}
$$
\begin{diagram}[w=72pt,h=60pt,tight]
  XIIY & \bigtopglob{X\lambda_I Y}
  \lift{15}{X\cel{K}^\lambda} 
  \straighthor{X\lambda_{IY}}
  \lift{-15}{\cel K\inv_{X,IY}}
  \bigdropglob{\rho_X IY}  & 
  XIY  \\
  \thinleftglob{X\lambda_I Y}\lift{0}{\cel{K}_{X,I}Y}\thinrightglob{\rho_X IY}  &  
  \comm
  & \leftglob{X\lambda_Y}\lift{0}{\cel K_{X,Y}}\rightglob{\rho_X Y}  \\
  XIY  & \rTo_{X\lambda_Y}  & XY
  \end{diagram}
  \vspace{24pt}
$$
Here finally the three $2$-cells incident to the $XIIY$ vertex cancel each other out, thanks
to Axiom~\eqref{TA2}, and in the end, remembering $\alpha=\lambda_I$, we get
  
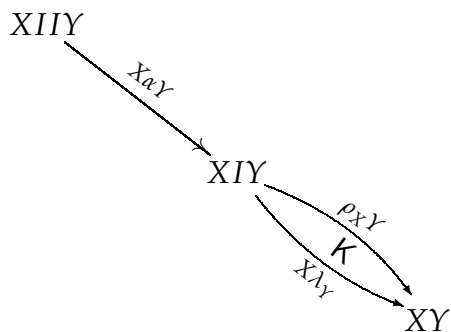
\vspace{24pt}
$$
\begin{diagram}[w=36pt,h=28pt,tight,hug]
  XIIY &&&&  \\
  &\rdTo^{X\alpha Y} &&&& \\
  && XIY && \\
  &&&\diagonaltopglob{\rho_X Y}
  \diagonaldropglob{X\lambda_Y}\NWSElabel{\cel K}   & \\
  &&&& XY
  \end{diagram}
  \vspace{24pt}
$$
as required to establish that $\cel K$ satisfies 
Equation~\eqref{RYXL}.  Hence we have proved (ii)\impl(i), and therefore 
altogether (ii)\impl(ii').  The converse, (ii')\impl(ii) follows now by 
left-right symmetry of the statements.  (But note that the proof via (i)
is not symmetric, since it relies on the definition $\alpha=\lambda_I$.
To spell out a proof of (ii')\impl(ii), use rather $\alpha=\rho_I$, observing
that the intermediate result (i) would refer to different 
 $\cel L$ and $\cel R$.)
\end{dem}

We have now given a construction in each direction, but both constructions 
involved choices.  With careful choices, applying one 
construction after the other in either way gets us back where we started.
It is clear that this should constitute an equivalence of $2$-categories.
However, the involved choices make it awkward to make the correspondence
functorial directly.  (In technical terms, the constructions are 
ana-$2$-functors.)
We circumvent this by introducing an intermediate $2$-category dominating 
both.  With this auxiliary $2$-category, the results we already proved
readily imply the equivalence.

\begin{blanko}{A correspondence of $2$-categories of units.}
  Let $\UU$ be following $2$-category.  Its objects are septuples
  $$
  (I,\alpha,\lambda,\rho,\cel L, \cel R,\cel K),
  $$
  with equi-arrows
  $$
  \alpha:II\to I, \quad 
  \lambda_X: IX \to X , \quad 
  \rho_X : XI \to X,
  $$
  (and accompanying naturality $2$-cell data),
  and natural invertible $2$-cells
  $$
  \cel L : I \lambda_X \Rightarrow \alpha X, \quad
  \cel R: X\alpha \Rightarrow \rho_X I, \quad
  \cel K: X\lambda_Y\Rightarrow \rho_X Y .
  $$
  These data
  are required to satisfy Equation~\eqref{RYXL} (compatibility of $\cel K$ with
  $\cel L$ and $\cel R$).  
%   (It follows from our results
%   that the induced auxiliary families $\cel K^\lambda$ and $\cel K^\rho$ will
%   satisfy Equation~\eqref{IIAB} as well as \eqref{TA2} and \eqref{TA3}.)

The arrows in $\UU$ from $(I,\alpha,\lambda,\rho,\cel L, \cel R,\cel K)$ to
$(J,\beta,\ell,r,\cel L', \cel R',\cel H)$
are quadruples
$$
(u, \cel{U}\loft, \cel{U}\rught, \cel U) ,
$$
where $u:I \to J$ is an arrow in $\CC$, $\cel{U}\loft$ and $\cel{U}\rught$
are as in \ref{GPS}, and $\cel U$ is a morphism of pseudo-idempotents from
$(I,\alpha)$ to $(J,\beta)$.  These data are required to satisfy 
Equation~\eqref{PK} (compatibility with Kelly cells)
as well as  Equations~\eqref{P} and \eqref{Q} in Lemma~\ref{unitmap} (compatibility
with the left and right $2$-cells).

Finally a $2$-cell from $(u, \cel{U}\loft, \cel{U}\rught, \cel U)$
to $(v, \cel{V}\loft, \cel{V}\rught, \cel V)$ is a $2$-cell
$$
\cel T : u \Rightarrow v
$$
required to be a $2$-morphism of pseudo-idempotents (compatibility with $\cel{U}$
and $\cel{V}$ as in \ref{semimonoid-transf}), and to satisfy Equation~\eqref{TXP}
(compatibility with $\cel{U}\loft$ and $\cel{V}\loft$)  as well as 
Equation~\eqref{TXQ} (compatibility with $\cel{U}\rught$
and $\cel{V}\rught$).

Let $\EE$ denote the $2$-category of cancellable-idempotent 
units introduced in \ref{main-def},
and let $\GG$ denote the $2$-category of GPS units of \ref{GPS}.
We have evident forgetful (strict) $2$-functors
\begin{diagram}[w=4.5ex,h=4.5ex,tight]
  && \UU \\
  & \ldTo^{\Phi} && \rdTo^{\Psi} \\
\EE &&&& \GG .
\end{diagram}
\end{blanko}

\stepcounter{satz}
\begin{satz}[Equivalence]\label{thmE}
  The $2$-functors $\Phi$ and $\Psi$ are $2$-equivalences.
  More precisely they are surjective on objects and strongly fully faithful
  (i.e.~isomorphisms on hom categories).
\end{satz}

\begin{dem}
  The $2$-functor $\Phi$ is surjective on objects by Lemma~\ref{LR} and
  Proposition~\ref{Kcoh}.  Given an arrow $(u,\cel{U})$ in $\EE$ and overlying
  objects in $\UU$, Lemma~\ref{unitmap} says there are unique $\cel{U}\loft$ and
  $\cel{U}\rught$, and Proposition~\ref{Kelly-compatibility} ensures the
  required compatibility with Kelly cells (Equation~\eqref{PK}).  Hence $\Phi$
  induces a bijection on objects in the hom categories.  Lemma~\ref{unit2map}
  says we also have a bijection on the level of $2$-cells, hence $\Phi$ is an
  isomorphism on hom categories.  On the other hand, $\Psi$ is surjective on
  objects by \ref{LRfromK} and Proposition~\ref{KLR}.  Given an arrow
  $(u,\cel{U}\loft,\cel{U}\rught)$ in $\GG$, Lemma~\ref{W} below says that for
  fixed overlying objects in $\UU$ there is a unique associated $\cel{U}$, hence
  $\Psi$ induces a bijection on objects in the hom categories.  Finally,
  Lemma~\ref{unit2map} gives also a bijection of
  $2$-cells, hence $\Psi$ is strongly fully faithful.
\end{dem}

\begin{lemma}\label{W}
  Given a morphism of GPS units
  \begin{diagram}[w=12ex,tight]
  (I,\lambda,\rho,\cel K) & \rTo^{(u,\cel{U}\loft,\cel{U}\rught)}  & 
  (J,\ell,r,\cel H)
  \end{diagram}
  fix an equi-arrow $\alpha: II \isopil I$ with natural families
  $\cel L_X: I\lambda_X \Rightarrow \alpha X$ and $\cel R_X: \alpha X 
  \Rightarrow \rho_X I$ satisfying Equation~\eqref{RYXL} (compatibility with 
  $\cel K$), and fix an equi-arrow $\beta: JJ \isopil J$ with natural families
  $\cel L'_X: I\ell_X \Rightarrow \beta X$ and $\cel R'_X: \beta X 
  \Rightarrow r_X I$ also satisfying Equation~\eqref{RYXL} (compatbility with 
  $\cel H$).  
  Then there is a unique $2$-cell
    \begin{diagram}
  II & \rTo^{uu}  & JJ  \\
  \dTo<\alpha  &  \cel{U}  & \dTo>\beta  \\
  I  & \rTo_u  & J .
  \end{diagram}
  satisfying Equations~\eqref{P} and \eqref{Q}
  (compatibility with $\cel{U}\loft$ and the
  left $2$-cells, as well as compatibility with $\cel{U}\rught$ and the right 
  $2$-cells). 
\end{lemma}

\begin{dem}
  Working first with left $2$-cells, define a family $\cel W_X$ by the 
equation
%%  generate horizontalCylinder \
%%     startobjects IIX IX \
%%     endobjects JJX JX \
%%     straightarrows {uuX} {uX} \
%%     curvedarrows {I\lambda_X} {\alpha X} {J\ell_X} {\beta X} \
%%     twocells {\cel L_X} {\cel L'_X} {\cel W_X} {u \cel{U}\loft_X} \
%%     label {}
%-------------------------------------output written 2009/07/14 20:15:00 CEST----
$$
\begin{diagram}[w=6ex,h=6ex,tight]
  IIX&\rTo^{uuX}&JJX\\
  \leftglob{I\lambda_X}\lift{0}{\cel L_X}\rightglob{\alpha X}&\shuft{20}{\cel W_X}&\rightglob{\beta X}\\
  IX&\rTo_{uX}&JX
\end{diagram}
\qquad\qquad = \qquad\qquad
\begin{diagram}[w=6ex,h=6ex,tight]
  IIX&\rTo^{uuX}&JJX\\
  \leftglob{I\lambda_X}&\shuft{-20}{u \cel{U}\loft_X}&\leftglob{J\ell_X}\lift{0}{\cel L'_X}\rightglob{\beta X}\\
  IX&\rTo_{uX}&JX
\end{diagram}
$$
%-------------------------------------end of output------checksum:0xBF09B993----
It follows readily from Equation~\eqref{KP=PH} that the family has the property
$$
\cel W_{XY} = \cel W_X Y
$$
for all $X, Y$, and it is a standard argument that since a unit object exists,
for example $(I,\lambda_I)$, this implies that
$$
\cel W_X = \cel{U} X
$$
for a unique $2$-cell 
  \begin{diagram}
  II & \rTo^{uu}  & JJ  \\
  \dTo<\alpha  &  \cel{U}  & \dTo>\beta  \\
  I  & \rTo_u  & J ,
  \end{diagram}
  and by construction this $2$-cell has the required compatibility with 
  $\cel{U}\loft$ and the left constraints.  
  To see that this $\cel U$ is also compatible with $\cel{U}\rught$ and the 
  right constraints we reason backwards: $(u,\cel U)$ is now a morphisms of
  units from $(I,\alpha)$ to $(J,\beta)$ to which we apply the right-hand version of
  Lemma~\ref{unitmap} to construct a new $\cel{U}\rught$, characterised by
  the compatibility condition.  By Proposition~\ref{Kelly-compatibility}
  this new $\cel{U}\rught$ is compatible with $\cel{U}\loft$ and the Kelly
  cells $\cel K$ and $\cel H$ (Equation~\eqref{PK}),
  and hence it must in fact be the original
  $\cel{U}\rught$ (remembering from \ref{GPS} that $\cel{U}\loft$ and
  $\cel{U}\rught$ determine each other via \eqref{PK}).  So the $2$-cell
  $\cel U$ does satisfy both the required compatibilities.
\end{dem}

\hyphenation{mathe-matisk}

% \bibliographystyle{scplain}
% \bibliography{joachims}

\begin{thebibliography}{10}

\bibitem{Benabou:bicat}
{\sc Jean B{\'e}nabou}.
\newblock {\em Introduction to bicategories}.
\newblock In {\em Reports of the Midwest Category Seminar}, no.~47 in Lecture
  Notes in Mathematics, pp. 1--77. Springer-Verlag, Berlin, 1967.

\bibitem{Gordon-Power-Street}
{\sc Robert Gordon, John Power, {\rm and }Ross Street}.
\newblock {\em Coherence for tricategories}.
\newblock Mem. Amer. Math. Soc. {\bf 117} (1995), vi+81pp.

\bibitem{Gurski:PhD}
{\sc Nick Gurski}.
\newblock {\em An algebraic theory of tricategories}.
\newblock PhD thesis, University of Chicago, 2006.
\newblock Available from http://www.math.yale.edu/\textasciitilde
  mg622/tricats.pdf.

\bibitem{Joyal-Kock:traintracks}
{\sc Andr{\'e} Joyal {\rm and }Joachim Kock}.
\newblock {\em Weak units and homotopy $3$-types}.
\newblock In A.~Davydov, M.~Batanin, M.~Johnson, S.~Lack, {\rm and }A.~Neeman,
  editors, {\em Ross Street Festschrift: Categories in algebra, geometry and
  mathematical physics}, vol. 431 of Contemp. Math., pp. 257--276. Amer. Math.
  Soc., Providence, RI, 2007.
\newblock ArXiv:math.CT/0602084.

\bibitem{Kelly:MacLanesCoherence}
{\sc Max Kelly}.
\newblock {\em On {M}ac{L}ane's conditions for coherence of natural
  associativities, commutativities, etc}.
\newblock J. Algebra {\bf 1} (1964), 397--402.

\bibitem{Kelly-Street:2cat}
{\sc Max Kelly {\rm and }Ross Street}.
\newblock {\em Review of the elements of {$2$}-categories}.
\newblock In {\em Category Seminar (Proc. Sem., Sydney, 1972/1973)}, no. 420 in
  Lecture Notes in Mathematics, pp. 75--103. Springer-Verlag, Berlin, 1974.

\bibitem{Kock:0507116}
{\sc Joachim Kock}.
\newblock {\em Weak identity arrows in higher categories}.
\newblock IMRP Internat. Math. Res. Papers {\bf 2006} (2006), 1--54.
\newblock ArXiv:math.CT/0507116.

\bibitem{Kock:0507349}
{\sc Joachim Kock}.
\newblock {\em Elementary remarks on units in monoidal categories}.
\newblock Math. Proc. Cambridge Philos. Soc. {\bf 144} (2008), 53--76.
\newblock ArXiv:math.CT/0507349.

\bibitem{MacLane:naturalAssociativity}
{\sc Saunders {Mac~Lane}}.
\newblock {\em Natural associativity and commutativity}.
\newblock Rice Univ. Studies {\bf 49} (1963), 28--46.

\bibitem{Saavedra}
{\sc Neantro Saavedra~Rivano}.
\newblock {\em Cat{\'e}gories Tannakiennes}.
\newblock No. 265 in Lecture Notes in Mathematics. Springer-Verlag, Berlin, 1972.

\bibitem{Simpson:9810}
{\sc Carlos Simpson}.
\newblock {\em Homotopy types of strict $3$-groupoids}.
\newblock ArXiv:math.CT/9810059.

\bibitem{Stasheff:1963}
{\sc Jim Stasheff}.
\newblock {\em Homotopy associativity of {$H$}-spaces. {I}, {II}}.
\newblock Trans. Amer. Math. Soc. {\bf 108} (1963), 275--292; 293--312.

\end{thebibliography}

\label{lastpage}

\end{document}